\documentclass[11pt]{article}   	                 % use "amsart" instead of "article" for AMSLaTeX format; if "amsart" 
\usepackage{amssymb, amsmath, amsthm}
\usepackage[colorlinks=true,linkcolor=blue,citecolor=red]{hyperref}
\usepackage{mathrsfs}                             % the \mathscr fonts
\usepackage[pdftex]{graphicx}                 % \includegraphics[scale =?, width=?\textwidth]{png/pdf/jpg/eps}
\usepackage[small]{caption}                     % Appearance of Caption, e.g. [font=small,labelfont=bf] 
\usepackage{wrapfig}                               %  \begin{wrapfigure}[lineheight]{position}{width}
\usepackage[subnum]{cases}
\usepackage{multirow}

\newtheorem{Lemma}{Lemma}[section]
\newtheorem{Theorem}{Theorem}
\newtheorem{Formal}{Formal Result}
\newtheorem{Proposition}[Lemma]{Proposition}

\newtheorem{Remark}[Lemma]{Remark}
\newtheorem{Definition}[Lemma]{Definition}

\newtheorem{Assumption}[Lemma]{Assumption}
\newtheorem{Example}[Lemma]{Example}

\newenvironment{Proof}[1][Proof]
 {\begin{trivlist} \item[]{\textbf{#1.} }}
{\hspace*{\fill}$\rule{.4\baselineskip}{.4\baselineskip}$\end{trivlist}}

\newenvironment{Acknowledgment}
 {\begin{trivlist}\item[]\textbf{Acknowledgments }}{\end{trivlist}}

\makeatletter\@addtoreset{figure}{section}\makeatother

\makeatletter \@addtoreset{equation}{section} \makeatother

\newcommand{\R}{\mathbb{R}}
\newcommand{\C}{\mathbb{C}}

\newcommand{\Z}{\mathbb{Z}}

\def\Re{\mathop{\mathrm{\,Re}\,}}

\def\diag{\mathop{\mathrm{\,diag}\,}}

\newcommand{\caA}{\mathcal{A}}
\newcommand{\caK}{\mathcal{K}}
\newcommand{\caL}{\mathcal{L}}

\newcommand{\caO}{\mathcal{O}}
\newcommand{\caE}{\mathcal{E}}

\newcommand{\rmd}{\mathrm{d}}
\newcommand{\rme}{\mathrm{e}}

\renewcommand{\span}{\mathrm{\,span}\,}

\renewcommand{\leq}{\leqslant}
\renewcommand{\geq}{\geqslant}

\newcommand{\rmnum}[1]{\romannumeral #1}
\newcommand{\Rmnum}[1]{\textrm{\uppercase\expandafter{\romannumeral #1\relax}}}

\def\beq{\begin{equation}}
\def\eeq{\end{equation}}
\def\veps{\vepsilon}
\def\veps{\varepsilon}
\def\u{\mathbf{u}}
\def\v{\mathbf{v}}
\def\w{\mathbf{w}}
\def\U{\mathbf{U}}
\def\H{\mathbf{H}}

\def\cad{\mathcal{D}}
\def\scb{\mathscr{B}}
\def\scn{\mathscr{N}}
\def\sca{\mathscr{A}}

\def\bm {\mathbf{m}}

\def \cF {\mathcal{F}}% bold-face in mathmode
\def\tP{\widetilde{P}}
\def \tW {\widetilde{W}}
\newfam\bifam
\font\tenbi=cmmib10 scaled \magstep1 \font\sevenbi=cmmib10 at 11pt
\font\fivebi=cmmib10 at 6pt \textfont\bifam = \tenbi
\scriptfont\bifam = \sevenbi \scriptscriptfont\bifam= \fivebi

\begin{document}

%%%%%%%%%%%%%%%%%%%%%%%%%%%%%%%%%%%%%%%%%%%%%%%%%%%%%%
% Title Page
%%%%%%%%%%%%%%%%%%%%%%%%%%%%%%%%%%%%%%%%%%%%%%%%%%%%%%
\title{Existence, Bifurcation, and Geometric Evolution of Quasi-Bilayers in the Multicomponent Functionalized Cahn-Hilliard Equation}
\author{Keith Promislow and Qiliang Wu\\
\textit{\small Department of Mathematics, Michigan State University}\\
\textit{\small 619 Red Cedar Road East Lansing, MI 48824}}

\date{}
\maketitle

\begin{abstract}
\noindent Multicomponent bilayer structures arise as the ubiquitous plasma membrane in cellular biology and as blends of 
amphiphilic copolymers used in electrolyte membranes, drug delivery, and emulsion stabilization within the context of synthetic chemistry. We present
the multicomponent functionalized Cahn-Hilliard (mFCH) free energy as a model which allows competition between bilayers with distinct composition and
between bilayers and higher codimensional structures, such as co-dimension two filaments and co-dimension three micelles.  We investigate the
stability and slow geometric evolution of multicomponent bilayer interfaces within the context of an $H^{-1}$ gradient flow of the mFCH, addressing the
impact of aspect ratio of the amphiphile (lipid or copolymer unit) on the intrinsic curvature and the codimensional bifurcation. In particular we derive a Canham-Helfrich sharp interface energy whose intrinsic curvature arises through a Melnikov parameter associated to amphiphile aspect ratio. 
We construct asymmetric homoclinic bilayer profiles via a billiard limit potential and show that co-dimensional bifurcation
is driven by the experimentally observed layer-by-layer pearling mechanism. 
%Under the $H^{-1}$ gradient flow of the mFCH free energy with the slow time scale $\tau =\veps^2t$, the set of 
% quasi-bilayer profiles are invariant and the evolution of quasi-bilayers reduces to an surface-area-preserving Willmore flow of the interface.
\end{abstract}

% \vfill

 \hrule
 {\small
 \begin{Acknowledgment}
 The first author acknowledges supported by the National Science Foundation through grant
 DMS-1409940. Both authors thank Brian Wetton for providing numerical simulations in support of the analysis 
 and Arjen Doelman for several beneficial discussions.
 \end{Acknowledgment}
% 
% {\bf Running head:} {Grain boundaries in the Swift-Hohenberg equation}
% 
% {\bf Corresponding author:} Arnd Scheel
% 
 {\bf Keywords:} %Multicomponent 
 Functionalized Cahn-Hilliard energy, Canham-Helfrich energy, Multi-component bilayer, Billiard limit. }

 %\newpage

%%%%%%%%%%%%%%%%%%%%%%%%%%%%%%%%%%%%%%%%%%%%%%%%%%%%%%
\section{Introduction}
%%%%%%%%%%%%%%%%%%%%%%%%%%%%%%%%%%%%%%%%%%%%%%%%%%%%%%
Amphiphilic molecules, such as lipids and functionalized polymers, are central to the self-assembly of intricately structured, 
solvent accessible nano-scaled morphologies and network structures. Indeed, Nafion, the generic ion separator in PEM fuel cells 
is comprised of a hydrophobic flourocarbon polymer-backbone functionalized by $SO_3H$ acidic side-chains. The result is a plastic-like ionomer 
that imbibes water to achieve solvent-accessible surface areas on the order of $1000 $m$^2/$gram, the highest for any manufactured material, \cite{KnoxVoth}. 
By combining hydrophilic and hydrophobic groups, this amphiphilic molecules, which we will refer to as amphiphiles, serve as surfactants, 
lowering interfacial energy, often to negative values, \cite{hayward_jacs08}, and self-assembling into a cornucopia of possible morphologies 
that balance packing entropy against hydrophilicity. The resulting materials have received extensive attention, for their wide  applications 
to pharmaceutics,  emulsion stabilization, detergent production and energy conversion devices 
\cite{amphiphilic2000, ameduri, heeger, diat, PromislowWetton_2009,batesjain_2004,  laschewsky,li_2004, lutz}. 
Of course amphiphiles play a central role in cellular biology, particularly in plasma membranes. These complex multicomponent systems 
are highly heterogeneous in the lipid distribution and composition, yielding clustering of particular lipids that is associated to intrinsic 
curvature of the membrane \cite{Brandizzi_2005, Laventis_2010, Sansom_2014}. Moreover, the interplay between species of lipids is also critical to the
structure of the endoplasmic reticulum \cite{Brandizzi},  membrane signaling and trafficking \cite{SimonsVaz_2004}, and has been proposed as a possible pathways to cell division in primitive single-celled organisms, \cite{Szostak}.  

\subsection{Multicomponent functionalized Cahn-Hilliard free energy}
\begin{figure}[!ht]
  \begin{center}
  \begin{tabular}{p{3.1in}}
  \includegraphics[width=0.3\textwidth]{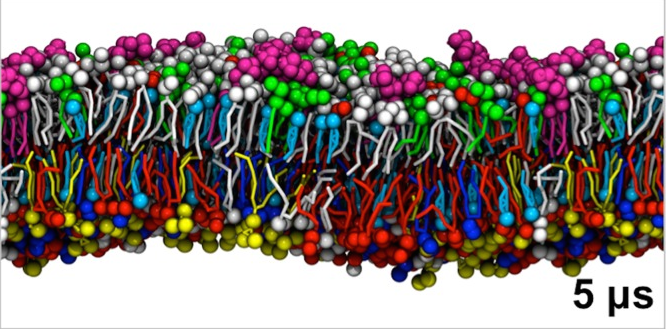}\\
   \includegraphics[width=0.3\textwidth]{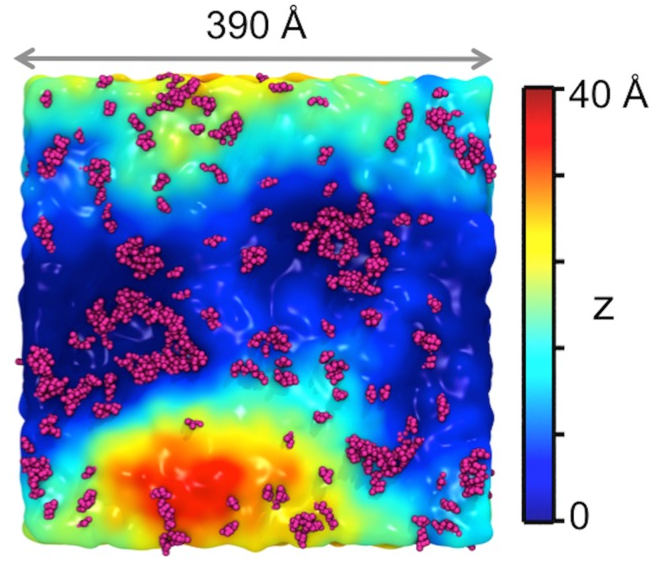}\\
   \includegraphics[width=0.3\textwidth]{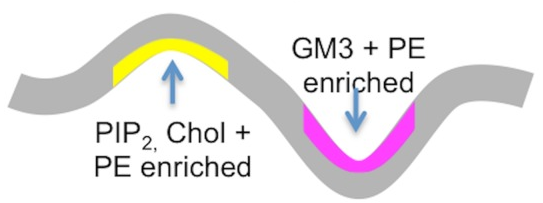}
  \end{tabular}
  % \vskip -2.2in
  \begin{tabular}{p{3.1in}}
  \includegraphics[width=0.44\textwidth]{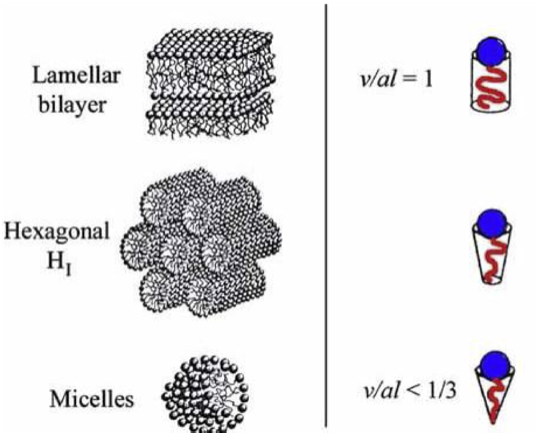} \\
    \includegraphics[width=0.44\textwidth]{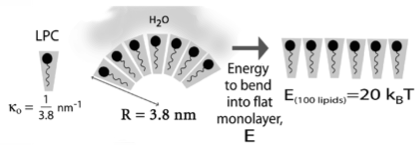} 
  \end{tabular} 
     \end{center}
     \vskip -0.3in
  \caption{(Left-panel) The initial and final structures of a plasma membrane from simulations, \cite{Sansom_2014}, composed of six types of lipids: POPC (light gray), POPE (red), POPS (blue), GM3 (magenta), Sph (green), PIP$_2$ (yellow), and Cholestoral (cyan), chemical descriptions are given in the appendix. The inner leaflet has 
initial composition POPC:POPE:Sph:GM3:Chol in ratios (40:10:15:10:25) while the outer leaflet had initial composition POPC:POPE:POPS:PIP$_2$:Chol in ratios (10:40:15:10:25).  The system was relaxed for $5\mu$ seconds, and the top image shows a snap-shot of the through-plane structure of the two leaflets, the middle image is a top down view color-coded according to the deformation of the membrane along the normal direction, with GM3 lipids depicted in magenta and other lipids not shown. 
 The bottom image relates the correlation between lipid type and membrane curvature. Reprinted from \cite{Sansom_2014}. PLoS Computational Biology is an open access journal.
 (Right-panel-top) Mouritsen's characterization of morphologies generated by lipid-solvent mixtures based upon the aspect ratio, $V/al$, of the truncated cone of volume $V$, cap area $a$, and length $l$, that best contains the hydrophobic (tail) component, \cite{Mouritsen}. 
 In the mFCH free energy, (\ref{e:mFCH1}), the aspect ratio is incorporated into both the intrinsic curvature, via the non-solenoidal perturbation ${\bf V}$, and the selection of co-dimension through pearling bifurcations mediated via the average value of $P(\u)$ over the profile,
 with co-dimension 3 micelles correspond to small values of $P$, and increasing values leading to co-dimension 2 filaments, and co-dimension 1 bilayers. Reprinted from \cite{Mouritsen} with permission from John Wiley and Sons.
 (Right-bottom) For LPC, \cite{FredCohen_2004}, the relation between the aspect ratio of the hydrophobic cone, the intrinsic curvature, and the energy required  to flatten the layer. 
 Reprinted from \cite{FredCohen_2004} with permission from Springer.
 } 
\label{f:packing}
\end{figure}

Many phase-field based models of cell membranes are based upon a ``single-layer'' formulation, \cite{DuWang, Lowengrub, Cohen_Ryham}, in which
a  co-dimension one interface separates two solvent phases, with the inside and the outside of the cell identified with a distinct phase field label. 
However many of the fundamental properties of membranes arise from their bilayer nature, in which an amphiphilic phase forms a two-sheeted co-dimension one 
interface that interpenetrates and separates a single solvent phase. Bilayer models of membranes have several advantages over single-layer models 
including a strong binding energy between the constituent layers which affords a natural mechanism to modulate the bilayer width, the possibility to 
perforate the membrane, or to reorganize into a higher co-dimensional  structure such as a filament or a micelle, see Figure \ref{f:packing}. 

We propose the multicomponent functionalized Cahn-Hilliard (mFCH) free energy as a model for the energy landscape of multicomponent blends of
amphiphiles (particularly lipids). The mFCH  free energy encompasses competition not only among morphologies with 
distinct co-dimensions but also couples the rearrangement of amphiphiles on the surface of network structure to the stability and geometric evolution of the
underlying surface. This work focuses on co-dimension one interfaces, and addresses  the competition with higher co-dimensional morphologies 
through the {\sl pearling} bifurcations that lead to filamentation or budding into micelles. The mFCH model is partially phenomenological, but incorporates the key 
elements of packing entropy and hydrophilicity in a transparent and mathematically elegant manner. Indeed, lipids have typically been classified based 
upon their aspect ratios, \cite{Mouritsen}, with bilayers associated with an aspect ratio near one, and smaller values of the aspect ratio associated with 
morphological bifurcations into filaments and micelles. However lipid aspect ratio also leads to non-zero intrinsic curvature of bilayers when the two sheets 
of the bilayer are formed of distinct mixtures. While the issues of intrinsic curvature and co-dimension bifurcation are both impacted by lipid aspect ratio, 
they have typically been studied independently, see \cite{Mouritsen} and \cite{Sansom_2014} for example. It is fundamental that a model distinguish between 
conditions under which an increase in aspect ratio of a lipid component would induce intrinsic curvature and when it would lead the bilayer to bifurcate into a 
higher co-dimensional structure, such as a filament or a micelle. Such a study requires a free energy which supports all these distinct phases, and we posit
that the multi-component functionalized Cahn-Hilliard (mFCH) is the simplest continuum model which does so.

The scalar functionalized Cahn-Hilliard (FCH) free energy is a well-studied model for mixtures of a single amphiphilic phase with a solvent, \cite{PromislowWetton_2009, GHYPromislow_2011, DaiPromislow_2013, PromislowYang_2014, doelmanpromislow_2014, PW_14}. The FCH corresponds to a special case of 
the model proposed by Tubner and Stray \cite{TS-87} and later Gommper and Shick \cite{GS-90};
motivated by small-angle X-ray scattering (SAXS) data of micro-emulsions of soapy-oil within water, they suggested a free energy for the soapy-oil volume fraction $u\in H^2(\Omega)$ 
in the form 
\beq
\cF_{\rm GS}(u) = \int_\Omega \veps^4\frac12  |\Delta u|^2 +\veps ^2 G_1(u) \Delta u + G_2(u)\, \rmd x, 
\eeq
where the function $G_1$ takes distinct values in the soapy-oil phase $u\approx 1$ and the water phase $u=0.$ The parameter $\veps>0$ scales homogeneously with space, and 
denotes the ratio of the interfacial thickness to the size of the domain $\Omega\subset \R^d$, generically $\veps\ll1$. This free energy includes a vast parameter space, and the functionalized Cahn-Hilliard corresponds to a simplifying choice $G_1=-W'(u)$, for $W$ a double-well potential with unequal depth wells at $u=0$ and $u=b>0$ and $G_2=\frac12 (W'(u))^2-\veps^p P(u)$. This form corresponds to a perturbation of a perfect square. Indeed, slightly redefining $W$ and $P$, see \cite{PW_14} for details, this special case of the Gompper and Schick free energy can be written as
the FCH free energy 
\beq
\label{e:FCH}
\cF(u)=\int_\Omega \frac{1}{2}\left(\veps^2\Delta u-W^\prime(u)\right)^2-\veps^p \left(\frac{1}{2}\veps^2\eta_1|\nabla u|^2+P(u)\right)\rmd x.
\eeq
The FCH admits a cornucopia of potential minimizers corresponding to bilayers, filaments, micelles and their local defects which, to leading order, render the dominant quadratic term zero;
that is they solve
 \beq
   \veps^2 \Delta u - W'(u) =O(\veps),
   \label{eq-CHcp}
  \eeq
 where the double-well  $W$ is the mixing potential, encoding the dominant components of the packing energy of the lipid-solvent mixture corresponding to {the} lipid volume 
 fraction $u$. These solutions are asymptotically close to the critical points of the corresponding Cahn-Hilliard type energy
 \beq
 \caE(u) = \int_\Omega \frac{\veps^2}{2} |\nabla u|^2 + W(u)\, \rmd x, 
 \eeq
 and correspond to optimal packings of lipid structures.  The functionalization terms are perturbations of the dominant quadratic, with the values $p=1$ and $p=2$
 corresponding to the strong and weak functionalizations respectively. In the strong functionalization these terms dominate the $O(\veps^2)$ Willmore residual in the 
 quadratic term, while for the weak functionalization the terms asymptotically balance. In both cases the functional terms encode perturbative energy preferences 
 among the optimal structures defined in (\ref{eq-CHcp}). For $\eta_1>0$ the first functionalization term rewards variation in $u$ associated with the generation of co-dimension one interface or higher codimension hyper-surfaces. This term encodes the strength of the hydrophilic-solvent interaction. The second functionalization term, given by the potential $P$, is typically taken in the simplified form $P=\eta_2 W$, incorporates energetic deviations among the solutions of (\ref{eq-CHcp}), assigning lower energies to morphologies residing in regions where $P>0.$
 
The scalar FCH free energy supports quasi-minimizing bilayer solutions that correspond to finite-width versions of the co-dimension-one sharp interface, $\Gamma$ immersed in $\Omega\subset R^d$.
In the $\veps\to0$ limit, the FCH evaluated at these bilayer solutions tends to a  Canham-Helfrich type limit \cite{canham_1970},
 \beq\label{e:canhel}
\caE_{\rm Helfrich}(\Gamma) := \int_\Gamma c_0(H(s)-c_1)^2 + c_2 + c_3 K(s)\, \rmd s
\eeq
where $H$ and $K$ are the mean and Gaussian curvatures of $\Gamma$, with $c_0,  c_1, c_2,$ and $c_3$ denoting constants; in particular $c_1$ is the intrinsic curvature, the value of the mean curvature that corresponds to zero bending energy, see \cite{GHYPromislow_2011} and Figure\,\ref{f:packing}. However for the scalar FCH, the bilayer must be symmetric about its center line and the intrinsic curvature, $c_1$, is zero in the sharp
interface limit. {Non-zero intrinsic curvature requires symmetry breaking}, and a bilayer that is asymmetric 
about its center line must be composed of more than one type of amphiphilic species. 

An  investigation of the interaction between intrinsic curvature in bilayers and morphological bifurcation, requires  a multi-component extension of the FCH free energy
that  encompasses a single solvent phase with $N$ lipid species. It is natural to introduce a vector-valued  phase function $u\in [H^2(\Omega)]^N$ with the $i$'th component denoting the local volume fraction of lipid type $i=1, \ldots, N$, and the associated multi-component  Gommper-Schick (mGS) free energy
\beq\label{e:mGS}
\cF_{mGS}(\u) = \int_\Omega \frac{\veps^4}{2} |{\bf D} \Delta \u|^2 + \veps^2 {\bf G}_1(\u)\cdot {\bf D} \Delta \u + G_2(\u) \rmd x,
\eeq
where the matrix, $\mathbf{D}=\diag(d_1, d_2, \dots, d_{N})>0$, encodes the differences in molecular length and weight between the
amphiphilic species. A key step in the selection of a minimal reduction of the mGS that preserves the richness of its solution space is the observation that 
for $N>1$ vector valued functions ${\bf G}: \R^N\mapsto\ R^N$ are not generically gradients of a scalar valued function. In the spirit of a minimal model, 
we take
 \beq\label{e:non-solenoidal}
 {\bf G}_1(\u):= -\nabla_{\u} W(\u) + \veps {\bf V}(\u), 
 \eeq
 for which the generic {\sl solenoidal} perturbation  ${\bf V}:\R^N\mapsto \R^N$ satisfies $\nabla_{\u}\times {\bf V}\neq 0$, in contrast to the irrotational term  $\nabla_{\u} W$ whose curl is zero. 
The most restrictive assumption on the form of the mGS, which parallels that made for the scalar case, is that the energy is close to a perfect square,
specifically, within the context of the weak functionalization, $G_2$ takes the form
\[
G_2(\u) := \frac12 \left| -\nabla_\u W (\u) +\veps {\bf V}(\u)\right|^2 - \veps^2 P(\u), 
\]
where the perturbation $P:\R^N\mapsto \R$.  Within this framework, the multi-component functionalized Cahn-Hilliard (mFCH) free energy
takes the form
\beq
\label{e:mFCH1}
\cF_{M}(\u)=\int_\Omega \frac{1}{2}|\veps^2{\bf D}\Delta \mathbf{u}-\nabla_\u W(\u)+\veps{\bf V}(\u)|^2-\veps^2 P(\u)\rmd x.
\eeq
To separate out the amphiphilicity term within the functionalization we make one last adjustment, redefining $W=\tW(\u)-\veps^2 A({\bf D}^{-\frac12}\u)$,
so that $\nabla_\u W = \nabla_u \tW-\veps^2 {\bf D}^{-\frac12}\nabla_\u A$, we may re-expand the quadratic term, which yields
\[
\cF_{M}(\u) = \int_\Omega \frac{1}{2}|\veps^2{\bf D}\Delta \mathbf{u}-\nabla_\u \tW(\u)+\veps{\bf V}(\u)|^2-\veps^2 \left(\tP(\u)-
\veps^2 {\bf D}^{\frac12}\Delta \u\cdot \nabla_\u A({\bf D}^{-\frac12}\u)\right) \rmd x,
\]
where $\tP=P+\nabla_\u \tW\cdot {\bf D}^{-\frac12}\nabla_\u A+\caO(\veps).$  Dropping the tilde notation, integration by parts on the last term yields the form
\beq
\label{e:mFCH3}
\cF_{M}(\u) = \int_\Omega \frac{1}{2}|\veps^2{\bf D}\Delta \mathbf{u}-\nabla_\u W(\u)+\veps{\bf V}(\u)|^2-\veps^2 \left(P(\u)+
\veps^2 \nabla \u :  {\bf D}^{\frac12}\nabla_\u^2 A {\bf D}^{-\frac12}\nabla \u\right) \rmd x,
\eeq
where the symbol $:$ denotes the double-contraction inner product that generates the usual norm on $\R^{N\times d}$.
In their density functional based model, \cite{andreussi_2012}, Andreussi, Dabo, and Marzari found that effective parameterization of 
solvation energies of solutes required two fit parameters associated to the quantum surface area and quantum volume associated to 
solvent-excluded region about the solute.  In this spirit, we reduce the functionalization terms  within the mFCH to an equivalent two-parameter family.
Generically we assume that $A$ is positive definite, so that its Hessian is a positive definite matrix, however for the two-parameter reduction
we further restrict $A$ to take the form
\[
A({\bf D}^{-\frac12}\u)= \frac12 \eta_1  |\u|^2,
\]
where the parameter  $\eta_1>0$ encodes the (common) strength of the hydrophilic interaction of the amphiphilic species with the solvent.

To better motive a simplified form for $P$, we will show in Section 2 that the co-dimension one bilayer morphologies are described to leading order
by solutions, $\phi_h$ of the second order system
\beq
\label{e:consprob}
 \partial_z^2\u = \nabla_{\bf u} {W}(\u),
 \eeq
 where $z$ is signed, scaled distance normal to the bilayer, and $u$ is  homoclinic to the solvent phase $\u=0.$ 
 The impact of $P$ is perturbative; for a quasi-critical point $\u_c$ of associated Cahn-Hilliard
free energy, that is a solution of (\ref{eq-CHcp}), the residual of the dominant quadratic term in (\ref{e:mFCH3}) arises at the same asymptotic order as the functionalization terms.
The volume integral of the funtionalized term,  $\veps^2\int_\Omega P(\u_c)\, \rmd x$, contributes at the same order as the Helfrich term in the quadratic residual, balancing packing entropy against geometry. The Hamiltonian structure of
(\ref{e:consprob}) requires that solutions homoclinic to zero reside on the set $\{\u \, \mid \, W(\u)\geq 0\}$ while higher co-dimensional profiles 
enter the region $\{ u \,\mid\, W(\u) <0\,\}$. This observation motivates the choice of $P$ as a scalar multiple of $W$, since volume integrals of $W(u)$ make a positive contribution for bilayers,
but give zero or negative contributions for codimension two and three morphologies.

\begin{Assumption}
\label{a:0} In the sequel we consider the mFCH, \eqref{e:mFCH3}, with the simplified form for the functionalization terms
\[
{\bf D}={\bf I}_N,\quad A({\bf D}^{-\frac12}\u)= \frac12 \eta_1 |\u|^2, \quad {\bf P}(\u) = \eta_2 W(\u), \quad {\bf V}(0)=0, \quad W(0)=0, 
\]
where $\eta_1\in\R^+,\eta_2\in\R$ are parameters. Moreover, the origin is a strict local minima of the mixing potential $W$, with a strictly positive-definite
Hessian:
$\nabla_u^2 W(0) > 0.$
\end{Assumption} 
Under Assumption \ref{a:0}, the weak multi-component weak FCH free energy takes the final form 
\beq
\label{e:mFCH}
\cF_{M}(\u) = \int_\Omega \frac{1}{2}|\veps^2\Delta \mathbf{u}-\nabla_\u W(\u)+\veps{\bf V}(\u)|^2-\veps^2 \left(
\veps^2 \frac{\eta_1}{2}| \nabla \u |^2  +\eta_2 W(\u) \right) \rmd x.
\eeq

The associated $H^{-1}$ gradient flow of the weak mFCH free energy is the system 
\beq
\label{e:mfch}
\u_t=\Delta \frac{\delta\cF_{M}}{\delta \u},
\eeq
subject to periodic boundary conditions on $\Omega\subset\R^d$, where the variational derivative takes the form
\begin{eqnarray}
 \frac{\delta\cF_{M}}{\delta \u} & = &\left(\veps^2\Delta-\nabla_\u^2W(\u)+\veps(\nabla_\u {\bf V}(\u))^{\rm T}\right)\left(\veps^2\Delta \u-\nabla_\u W(\u)+\veps {\bf V}(\u)\right) \nonumber \\
                                                             && +\,\veps^2\left(\veps^2\eta_1\Delta \u-\eta_2\nabla_\u W(\u)\right).
 \end{eqnarray}
The flow is local in space and conserves the $\veps$-scaled total mass-vector 
\beq
\mathcal{M}(\u):=\frac{1}{\veps}\int_\Omega \u\rmd x,
\label{e:calM-def}
\eeq
of the $N$ amphiphilic species.

%%%%%%%%%%%%%%%%%%%%%%
\subsection{Quasi-bilayers}
%%%%%%%%%%%%%%%%%%%%%%
Our analysis addresses the construction, linear stability, and slow evolution of families of \textit{quasi-minimizers}:  distributions $\u\in [H^2(\Omega)]^N$ with sufficiently small free energy
and low lipid mass fraction, $\mathcal{M}$ see (\ref{e:calM-def}). More specifically,  given a fixed constant $C>0$, the associated set of quasi-minimizers is given as 
\beq
\label{e:quamin}
\mathcal{Q}_C:=\left\{\u\in [H^2(\Omega)]^N \, \Bigg| \, |\mathcal{M}(\u)|\leq C \text{ and } \mathcal{F}_{M}(\u)\leq C\veps^3\right\}.
\eeq
The energy $\cF_{M}$ is bounded from below for fixed values of $\veps>0$, but can be negative, see \cite{keithzhang_2013} %Paper [52] on my web page
for a discussion of lower bounds for the scalar FCH.

Candidates for quasi-minimizers are readily constructed from approximate solutions of the stationary weak mFCH equation,
\beq
\label{e:smfch}
\frac{\delta \cF_M}{\delta \u}=\veps^2\bm,
\eeq
where $\bm\in \R^N$ can be viewed as an $\veps^2$-scaled Lagrange multiplier associated to the mass constraint imposed by the mFCH gradient flow. We focus our attention on
a class of quasi-minimizers called \textit{quasi-bilayers}. These are constructed by fixing an admissible co-dimension one base 
interface $\Gamma$ immersed in $\Omega$, changing to local coordinates $(z,s)$, where $z$ is $\veps$ scaled distance to $\Gamma$ and $s:S\mapsto \Gamma$ is a parameterization of $\Gamma$, see Definition \ref{d:noselfin} for details. Within the local variables the Laplacian takes the form (\ref{e:lap}), and the stationary equation (\ref{e:smfch}) reduces at lowest order to the second order dynamical system
 \beq
 \label{e:leadode}
 \partial_z^2\phi_h-\nabla_\u W(\phi_h)=0,
 \eeq
subject to the condition that $\phi_h$ is homoclinic to zero. For special cases, the homoclinic solution can be corrected to yield a higher-order approximation to (\ref{e:smfch}), however
the solenoidal perturbation, $\veps V$ is not such a case, and persistence of homoclinic solutions under the solenoidal perturbation requires the
introduction of a Melnikov parameter. To account for the necessary degrees of freedom within the system, we introduce the perturbed 
homoclinic $\Phi_h=\Phi_h(z,\veps; \bm)$, of the  solenoidal ODE,
\beq
\label{e:sole-ode}
\partial_z^2 \Phi_h+\veps a\partial_z\Phi_h -\nabla_\u W(\Phi_h)+\veps {\bf V}(\Phi_h)={-\veps^2 \left[\nabla_\u^2 W(0)\right]^{-1} \bm},
\eeq
where the mass constraint $\bm$ is viewed as a parameter that prescribes the far-field behavior
\beq
\Phi_\infty(\veps;{\bf m}):=\lim\limits_{z\to\infty} \Phi_h(z,\veps;\bm) = \veps^2[\nabla_\u^2W(0)]^{-2}\bm + O(\veps^3).
\eeq
Conversely the Melnikov parameter $a$ is fixed by the choice of the homoclinic profile $\phi_h$ and more significantly by $V$ through the expansion 
\beq
\label{e:aexp}
a=a_0+\caO(\veps), \quad a_0:=-M_1^{-1}\int_\R  {\bf V}(\phi_h(z))\cdot \partial_z\phi_h(z) \rmd z, \quad M_1:=\int_\R |\partial_z\phi_h(z)|^2\rmd z.
\eeq 
The \textit{quasi-bilayer} profiles are the vector-valued smooth functions $\u_q\in [H^2(\Omega)]^N$ obtained by the ``dressing'' of an admissible interface 
$\Gamma$, defined in Definition \ref{d:noselfin}, with the profile $\Phi_h(z,\veps;{\bf m})$ via the relations
\beq\label{e:quasib}
\u_q(x, \veps; {\bf m},\Gamma)=\begin{cases}
\Phi_h(z(x),\veps;{\bf m}), & x\in\Gamma_{l_0},\\
(1-\chi(\frac{\veps |z(x)|}{l_0}))\Phi_\infty(\veps; {\bf m})+\chi(\frac{\veps |z(x)|}{l_0})\Phi_{h,\delta}(z(x),\veps;{\bf m}), & x\in \Gamma_{3l_0}\backslash \Gamma_{l_0},\\
\Phi_\infty(\veps;{\bf m}), & x\in\Omega\backslash \Gamma_{3l_0},
\end{cases}
\eeq
where the inner region $\Gamma_{l_0}$ is given in Definition \ref{d:noselfin}, and $\chi:\R\to\R$ is a smooth cut-off function satisfying 
$\chi(r)=1$ for $|r|\leq 1$ and $\chi(r)=0$ for $|r|\geq 3$.  
%$$\Phi_\infty(\veps;{\bf m}):=\lim_{z\to\infty}\Phi_h(z,\veps;{\bf m})=\veps^2 [\nabla_\u^2W(0)]^{-2}{\bf m}+\caO(\veps^3).$$
The quasi-bilayer profiles are parameterized by $\{{\bf m}, \Gamma\}$, where ${\bf m}$ controls the far-field value and $\Gamma$ is an admissible co-dimension 
one immersion within $\Omega$ that specifies the bilayer center-line. For fixed $C, l_0>0$ we study the set of quasi-bilayers $\mathscr{M}_{C, l_0}(\veps)$ defined as 
\beq\label{e:setquasi}
\mathscr{M}_{C,l_0}(\veps):=\left\{ \u_q(\cdot, \veps; {\bf m}, \Gamma)\in C^\infty(\Omega) \mid \, |{\bf m}|\leq C,\, \Gamma\, \text{is an admissible interface with reach}\, l_0  \right\}.
\eeq

%%%%%%%%%%%%%%%%%%%%%%
\subsection{Main results}
%%%%%%%%%%%%%%%%%%%%%%

In Section \ref{s:2}, we
establish the persistence of homoclinic solutions under the solenoidal perturbations and verify that the set $\mathscr{M}_{C,l_0}(\veps)$ of quasi-bilayers
are in fact quasi-minimizers, and moreover their free energy $\cF_M(\u_{q})$ admits a Canham-Helfrich sharp interface free energy at leading order. 
More specifically, introducing the notation 
\beq\label{e:not}
{\bf M} := \int_\R\phi_{h}(z)\rmd z, \quad {\bf B}:=[\nabla_\u^2 W(0)]^{-2}{\bf m}, \quad \phi_{h,1}:=\partial_\veps \Phi_h|_{\veps=0},
\eeq
which denote, respectively the mass per unit length of the bilayer $\phi_h$, the leading-order constant far-field value of $\Phi_h$, and the $O(\veps)$-order term 
in the expansion of $\Phi_h$, we establish the following Theorem.
\begin{Theorem}\label{t:canhel}
Under Assumption \ref{a:homoex}, for any $C, l_0>0$ there exists $\veps_0, C_1>0$ such that for all $\veps\in(0,\veps_0)$ the family of quasi-bilayers $\mathscr{M}_{C,l_0}$ 
is contained within the set of quasi-minimizers  $\mathcal{Q}_{C_1}$ and the mFCH evaluated at $\u_q({\bf m}, \Gamma)\in\mathscr{M}_{C,l_0}$ reduces to a Canham-Helfrich sharp
interface energy on $\Gamma$ of the form
\beq\label{e:quasiFCH}
\mathcal{F}_M(\u_q)=\veps^3\frac{M_1}{2}\int_S\left[(H_0(s)-a_0)^2-(\eta_1+\eta_2)\right]\rmd s+\caO(\veps^4),
\eeq
where $H_0$ is the mean curvature of  $\Gamma$ and the intrinsic curvature $a_0$ and constant $M_1$ are defined in \eqref{e:aexp}. Moreover the total $\veps$-scaled mass takes the form
\beq\label{e:masscon}
\mathcal{M}(\u_q)= |\Gamma|{\bf M}+\veps\left[|\Omega|{\bf B}+|\Gamma|\int_\R\phi_{h,1}(z)\rmd z+\int_SH_0(s)\rmd s\int_\R z\phi_h(z)\rmd z\right]+\caO(\veps^2),
\eeq
\end{Theorem}
where ${\bf M}$, ${\bf B}$ and $\phi_{h,1}$ are defined in \eqref{e:not}.
%fixing the total mass  $\mathcal{M}$ of each amphiphile, the background state can be tuned by varying the base interface, or vice versa. 
%Numerical simulation shows, the condition that the total mass is far from being proportional to $\int_\R\phi_{h}(z)\rmd z$ drives the morphology into transient pearling right away. Figures.

%\vskip -2in
\begin{figure}[!ht]
%\vskip -0.3in
  \begin{center}
    \includegraphics[width=1\textwidth]{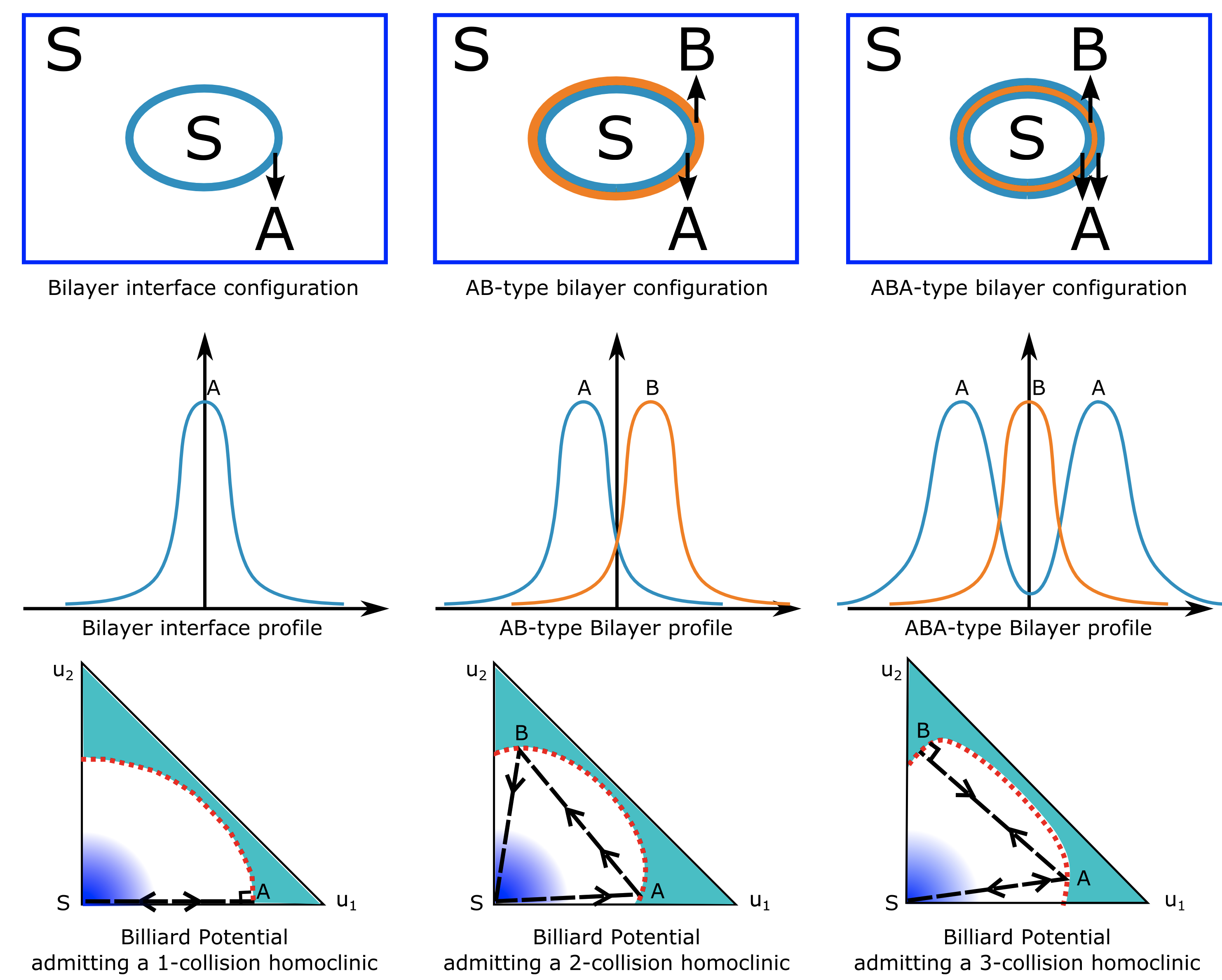}
     \end{center}
     \vskip -0.26in
  \caption{The top row depicts the physical domain $\Omega\subset\R^3$ with a bilayer (left), an AB-type bilayer (middle), and an ABA-type bilayer (right). 
  The solvent phase is denoted by S and A and B denote species of amphiphiles. The middle row depicts the density profile of the solute or solutes (phase A and B) 
 along the scaled normal direction, $z$, of the interface given by solutions of (\ref{e:consprob}). The bottom row depicts the mixing potential $W$ as a function of $\u$, with the superimposed dashed lines depicting the homoclinic orbits, $\phi_h$. The potentials are Birkhoff-billiard potentials,  see Definition \ref{d:bill} and Figure \ref{f:billiard} in Section 3 for more details.}
\vskip -0.14in
\label{f:SBlayer}
\end{figure}
%\vskip -2in

\begin{figure}[!ht]
\begin{center}
 % \begin{minipage}{0.2\textwidth}
  \begin{tabular}{lp{0.4in}}
  \includegraphics[width=0.242\textwidth]{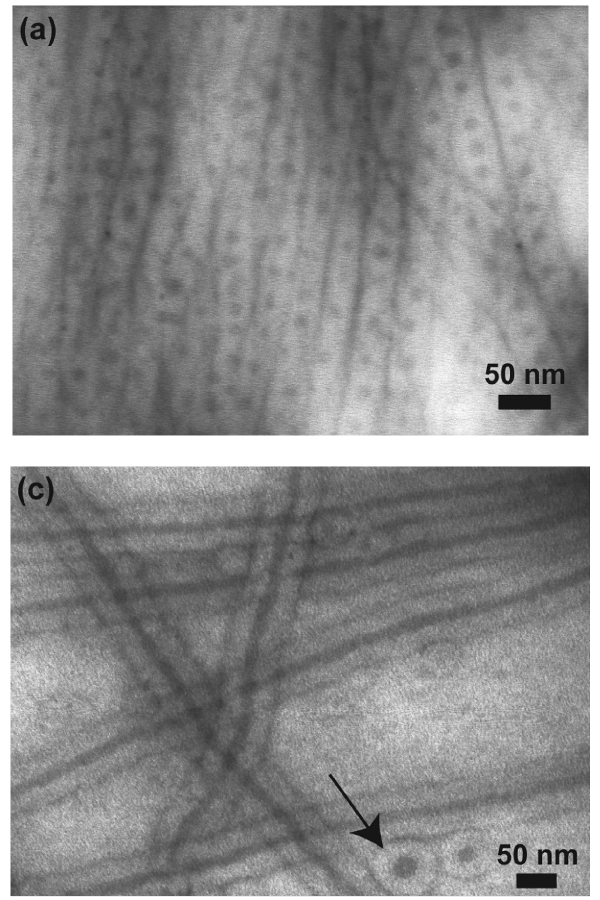} &
  \end{tabular}
%  \end{minipage}
%  \begin{minipage}{0.78\textwidth}
  \begin{tabular}{cc}
  \includegraphics[width=0.24\textwidth]{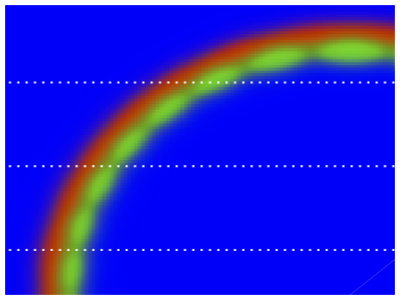} &  \includegraphics[width=0.18\textwidth]{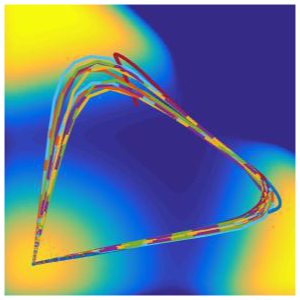} \\
 \includegraphics[width=0.24\textwidth]{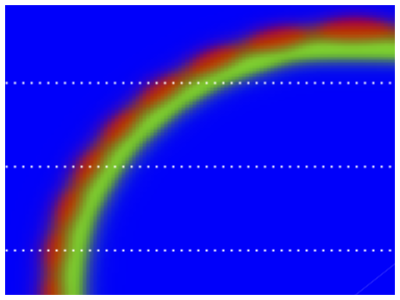} & \includegraphics[width=0.18\textwidth]{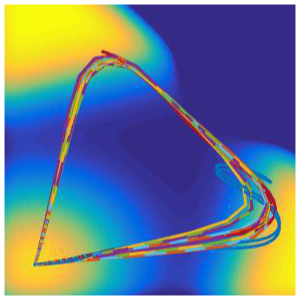} 
  \end{tabular}
 % \end{minipage}
 \end{center}
 %    \vskip -0.3in
  \caption{(Left)Distinct $AB$-type bi-layer filaments, formed from blends of PS$_{9.5K}$-PEO$_{5K}$ with PS$_{56K}$-P2VP$_{21K}$, exhibit a layer-by-layer pearling bifurcation into pearled P2VP cores. Reprinted with permission from \cite{hayward}. Copyright 2008 American Chemical Society. (Center, Right) Numerical simulations of the layer-by-layer pearling of an $AB$-type bilayer structure.  More specifically, given $N=d=2$, $\Omega=[0,2\pi]\times[0,2\pi]$, $\veps=0.2$, $\eta_1=\eta_2=1$, $V(\u)=(-u_2,u_1)^{\rm T}$ 
%for the first simulation and $V(\u)=(u_2,-u_1)$ for the second simulation 
and the mixing well as a regularized-billiard potential with $\delta=0.2$, the simulation of the weak mFCH gradient flow 
 with the initial data that is slightly above, equal to, and slightly below, the target value. The left panels show the physical configurations of the pearled morphologies, 
 where blue, red and green represent respectively the solvent, amphiphile A and amphiphile B. The right panels depicts the orbit of $\u$ parametrized by $x_1$ for chosen fixed $x_2$ (white dotted lines). The code used in the numerical simulation is developed by Brian Wetton. }
\label{f:pearling}
\end{figure}

Section \ref{s:3} focuses on a class of regularized \textit{Birkhoff billiard} potentials for the mixing wells for which it is possible to construct and explicitly characterize a large
family of homoclinic solutions of the corresponding Hamiltonian ODE \eqref{e:leadode}.  The regularized Birkhoff-billiard potentials $\{B_\delta(\u)\}$, are parameterized by a 
free parameter $\delta>0$, which is small in comparison to the features of $W$, and which controls the smoothing of the discontinuous Birkhoff potential. The Birkhoff billiard potential
$B(\u)$ is piece-wise smooth, with a strict local minima at the origin and transitions to a positive constant in a large region, and jumps discontinuously to a fixed negative 
value across a `collision-curve', see the bottom row of images of Figure\,\ref{f:SBlayer}. The Hamiltonian structure of (\ref{e:consprob}) precludes the orbits homoclinic to the origin from entering the negative-$W$ region and thus they impact the discontinuity in a billiard-like collision.  To each collision we associate a ``striation'' within the corresponding 
co-dimension one interface, see Figure~\ref{f:SBlayer}.
%the interactions between types of lipids lead to bilayer morphologies with hydrophilic head groups exposed to solvent and hydrophobic tail
%groups residing along the center line of the bilayer, \cite{batesjain_2004} [not a very good citation for this/ should we move this closer to Fig 1.4?], 
Striations have been observed experimentally in blends of polymer blends of $A-B$ diblock with a $B-C$ diblock, \cite{hayward}, with $A$ representing the most 
hydrophilic component and $C$ the most hydrophobic. These multicomponent, di-block blends phase-separate into distinct striated structures which exhibit striation localized pearling bifurcations, as depicted in the left panel of Figure \ref{f:pearling}.

Within the context of a regularized-billiard potential  $B_\delta(\u)$, we show the existence of $n$-striation bilayers which are homoclinic to the origin, 
and characterize the spectrum of the associated second variational derivative of $\cF_M$. We emphasize that the term bilayer refers to the dressing
of a co-dimension one admissible interface with a solution $\Phi_h$ of (\ref{e:sole-ode}) that is asymptotically homoclinic to the origin. The $n$ in the $n$-striation
refers to the number of distinct shifts in composition of the bilayer in the through-plane direction. In particular the term striation is only meaningful when the mixing
potential $W$ is a $\delta\ll1$ regularized Birkhoff billiard potential.  For a given regularized Birkhoff-billiard potential, the second variational derivative of the 
mFCH free energy  \eqref{e:mFCH} about an $n$-striation bilayer $\psi_{h,\delta}$ yields an operator $\mathbb{L}$ of the form
\beq\label{e:bbL}
\mathbb{L}:=\frac{\delta^2 \cF_M}{\delta \u^2}(\psi_{h,\delta})=(\caL_\delta+\veps^2\Delta_s)^2+\caO(\veps),
\eeq
where $\Delta_s$ is the Laplace-Beltrami operator associated to $\Gamma$ and the structural operator associated to $\psi_{h,\delta}$
%\beq
%\begin{matrix}
%\caL_\delta:& [H^2(\R)]^2 & \longrightarrow & [L^2(\R)]^2 \\
%& \u & \longmapsto & \partial_z^2\u-\nabla_\u^2 B_\delta(\psi_{h,\delta})\u.
%\end{matrix}
%\label{e:L-delta}
%\eeq
\beq
\caL_\delta:=  \partial_z^2 - \nabla_\u^2 B_\delta(\psi_{h,\delta}),
\label{e:L-delta}
\eeq
maps $[H^2(\R)]^2$ to  $[L^2(\R)]^2$. As for the scalar FCH in the weak \cite{hayprom_2014} settings, the negative spectrum of 
$\mathbb{L}$ can only arise when the dominant quadratic term is close to zero. Since the Laplace-Beltrami operator is non-positive, 
negative spectrum occurs only for tensor product eigenfunctions with the leading order form
\beq
\Psi_{j,k}(z,s):= \psi_j(z)\Theta_k(s) +O(\veps),
\eeq
where $\psi_j$ is an eigenfunction of $\caL_\delta$ associated to eigenvalue $\lambda_j\geq 0$ and $\theta_k$ is a Laplace-Beltrami eigenfunction
with eigenvalue $\beta_k<0$ for which $\lambda_j+\veps^2\beta_k=O(\sqrt{\veps}).$ Rigorous statements about the pearling spectrum
for the second variation, $\mathbb{L}$, for the scalar FCH free energy can be found in \cite{hayprom_2014}.

These tensor-product eigenmodes of $\mathbb{L}$ are the pearling
eigenvalues, as their high-frequency in-plane oscillation generates a pearled morphology, see Figure\,\ref{f:pearling} for experimental observations and 
\cite{PW_14} for a rigorous construction of pearled morphologies for the scalar FCH. 
%the nonnegative linear spectra associated to a quasi-bilayer is controlled by the nonpositive spectrum of $\mathbb{L}$. 
%Actually, Lemma 4.1 and 4.2 in \cite{doelmanpromislow_2014} can be verified for the multi-component case without much difficulty. 
%In addition, because of the positive nature of the dominant quadratic term in \eqref{e:bbL}, it follows that negative spectra of $\mathbb{L}$ can only be $\caO(\veps)$, and from the negative nature of $\Delta_s$, they must be associated to positive spectra of $\caL_\delta$.
For $\delta$ sufficiently small we show that $\caL_\delta$ associated to an $n$-collision homoclinic possesses precisely $n$ large positive eigenvalues of 
order $\caO(\delta^{-2})$ which may be enumerated so that for $j=1, \ldots, n$, the corresponding $j$'th eigenfunction $\{\psi_j\}_{j=1}^n$ is generically 
localized on the $j$-th striation, as corresponds to the experimentally observed striation-localized pearling. More precisely, we establish the following theorem.
%------------------------------------------------------------
\begin{Theorem}\label{t:eigen}
Given a Birkhoff-billiard potential $B(\u)$ in the form Definition \ref{d:bill} for which the leading-order ODE \eqref{e:leadode} admits an $n$-collision homoclinic $\psi_h$ satisfying Assumption 
\ref{a:homo} and \ref{a:trans}. For sufficiently small $\delta>0$, the leading-order ODE \eqref{e:leadode} with the regularized Birkhoff-billiard potential $B_\delta(\u)$ given in Definition \ref{d:quasi} admits a unique homoclinic continuation, $\psi_{h,\delta}$.  The associated linearized operator $\caL_\delta$ defined in (\ref{e:L-delta}) has 
real spectrum $\sigma(\caL_\delta)\subset\R$, and admits precisely $n$ ``collision'' eigenvalues, $\{\lambda_1(\delta), \cdots, \lambda_n(\delta)\}$, of order $\caO(\delta^{-2})$. 
More precisely, there exists $C_0>0$, $\delta_0>0$ and $c(\delta):(0,\delta_0]\to (0, C_0]$ with $\lim_{\delta\to 0+}c(\delta)=0$ so that for any $\delta\in(0,\delta_0]$ and $c\in(c(\delta), C_0]$,
\[
\sigma(\caL_\delta)\cap\{\lambda\in\C\mid \Re\lambda>c\delta^{-2}\}=\{\lambda_1(\delta),\lambda_2(\delta),\cdots,\lambda_n(\delta)\}.
\]
Moreover the collision eigenvalues may be labeled so that $\lambda_j(\delta)$ admits the expansion
\[
\lambda_j=\nu_j\delta^{-2}+\caO(\delta^{-1}),
\]
where $\nu_j$ is the unique positive eigenvalue of the $j$-th collision operator $\caK_j$ as defined in \eqref{d:limcol}. For the generic case when $\nu_i\neq \nu_j$ for $i\neq j$, then $\caL_\delta$ has $n$ distinctive collision eigenvalues, each of which is simple with the corresponding $j$'th eigenfunction $\psi_j$ localized on the $j$-th collision interval
$K_j(\delta)$ as defined in \eqref{e:cin}. In addition, we have
\[
\sigma_{ess}(\caL_\delta)=(-\infty,-2], \quad 0\in\sigma_{pt}(\caL_\delta),
\]
where, by Assumption \ref{a:trans}, $\lambda=0$ is simple with eigenspace spanned by the translational mode $\partial_z\phi_{h,\delta}$.
\end{Theorem}

In Section\,\ref{s:4} we return to the general mixing well, $W$, and use a formal matched asymptotic expansion to show that 
the manifold of quasi-bilayers is approximately forward invariant on the $O(\veps^{-2})$ time-scales, and derive the evolution of the underlying interface $\Gamma$
and background state $\bf{m}.$  
%This evolution corresponds to the reduction of the gradient flow (\ref{e:mfch}) onto the 
%center-stable manifold corresponding to the translational eigenvalues of the second variational derivative of  $\cF_M$ arising from the tensor product of the translational eigenfunction in the %kernel of $\caL$ with the high-frequency Laplace-Beltrami eigenmodes  of $\Gamma$.  
Compared with the scalar case \cite{DaiPromislow_2013}, the novelty of this derivation lies in the following two points. First, the mFCH equation preserves total mass, but the dressing of
an interface $\Gamma$ with a bilayer profile specifies the total mass, to leading order, as a multiple of the vector-valued mass per unit area ${\bf M}(\phi_h)\in\R^N$,  see (\ref{e:not}).
Consequently, generic initial data cannot converge to a quasi-bilayer profile if the mass of the initial data is not properly tuned. This effect is manifest in Lemma 4.1, 
on the $O(\veps^{-1})$ time scale for initial data which have an $O(\veps)$ mis-fit with the quasi-bilayer equilibrium -- the background state ${\bf m}$ starts at $O(\veps)$ but 
converges to $O(\veps^2)$ on the faster time scale only if the total mass is properly tuned. An analysis of pearling bifurcations for the scalar case shows a sensitivity to the value of the
far-field equilbrium, and this is manifested in our numerical investigations, depicted in the right-two panels of Figure\,\ref{f:pearling}. For a mixture of two lipid types and a solvent,
initial data that is off by $O(\veps)$ from the target ratio corresponding to the bilayer leads to a pearling bifurcation in the inner (respectively outer) layer of the bilayer depending upon the 
nature of the mis-match between initial and target mass vectors.  Lemma 4.1 shows that for initial composition corresponding to a 
multiple of the bilayer mass vector, the far-field equilibrium ${\bf m}$ tends to $O(\veps^2)$ on the fast $\caO(\veps^{-1})$ time-scale. Subsequently
we enter a slow evolution on the $O(\veps^{-2})$ time-scale corresponding to an area-preserving Willmore geometric flow, with the far-field state 
algebraically slaved to enforce conservation of total mass of each species.  The second novelty resides in the appearance of the intrinsic curvature in the 
Willmore formulation of the slow-time geometric flow.  More specifically, noting that the generalized mean curvature of the interface admits the expansion 
$H(z,s,\tau;\veps)=H_0+\veps H_1+\caO(\veps^2)$ and introducing respectively the squared bilayer mass and the curvature-weighted projection
\beq\label{e:not2}
M_2:=\int_\R|\phi_h(z)|^2\rmd z,  \quad \Pi_\Gamma(f):=f-H_0\frac{\int_S fH_0\rmd s}{\int_S H_0^2\rmd s},
\eeq
we use the method of multiple scales to establish the following formal result.
\begin{Formal}\label{t:willmore}
For small $\veps>0$, the method of matched asymptotics applied to  the $H^{-1}$ gradient flow of the mFCH free energy shows that in the slow time scale $\tau_2=\veps^2t$, 
a quasi-bilayer $\u_q$ associated to $\{{\bf m}_0, \Gamma_0\}$ evolves as a quasi-bilayer plus a small perturbation, that is,
\beq
\u(\tau_2,x,\veps)=\u_q(x,\veps;{\bf m}(\tau_2), \Gamma(\tau_2))+\u_p(\tau_2,x,\veps),
\eeq
where $\|\u_p\|_{L^\infty(\Omega)}=\caO(\veps^2)$ and $\|\u_p\|_{L^1(\Omega)}=\caO(\veps^3)$. More importantly, 
the evolution reduces to a surface-area-preserving Willmore flow of the interface $\Gamma(\tau_2)$, that is, for $\tau_2\geq0$, the total interfacial surface area
$|\Gamma(\tau_2)|$ equals the initial surface area $|\Gamma_0|$ to leading order, and the leading order normal velocity $V_{n,0}$ satisfies
\begin{equation}\label{e:willmore}
V_{n,0}(s,\tau_2)=M_2^{-1}M_1\Pi_\Gamma\left[  \left(\Delta_s-\frac{1}{2}H_0(H_0-a_0)-H_1  \right)  (H_0-a_0)\right],
\end{equation}
while the background state ${\bf m}(\tau_2)$ is slave to the Willmore flow via mass conservation 
\begin{equation}
\mathbf{m}(\tau_2)\cdot \mathbf{M}=-M_1\frac{\int_S\left[ -|\nabla_s H_0|^2 +\frac{\eta_1+\eta_2}{2} H_0^2-\frac12 H_0^2(H_0-a_0)^2-H_1H_0(H_0-a_0) \right] \rmd s}{\int_S H_0^2\rmd s},
\end{equation}
where
% the variables $V_{n,0}$ and $\gamma_0$ are respectively the leading-order term of the interface's normal velocity and area with respect to $\veps$, 
the Melnikov parameter $a_0$ and the bilayer parameter $M_1$ are defined in \eqref{e:aexp}, the mass per unit length ${\bf M}$ is defined in \eqref{e:not} and 
the curvature-weighted projection $\Pi_\Gamma$ and squared bilayer mass $M_2$ are defined in \eqref{e:not2}, .
%In addition, 
%for $\tau\geq 0$, the co-dimensional one morphology admits the expansion
%\[
%\u(x,\tau;\veps)=
%\begin{cases}
%\phi_h(z(x))-\veps\tilde{\u}_1(z(x))+\veps^2\tilde{\u}_2(x,\tau)+\caO(\veps^3),& x\in\Gamma_{l_0}(\tau),\\
%\veps^2\left[\nabla_\u^2W(0)\right]^{-2}{\bf B}_2(\tau)+\caO(\veps^3), &x\in \Sigma\backslash\Gamma_{l_0}(\tau),\\
%\end{cases}
%\]
%where $\tilde{\u}_1:=\caL_{0,\perp}^{-1}(a_0\partial_z\phi_h+{\bf V}(\phi_h))$, $\tilde{\u}_2\perp \partial_r\phi_h$ satisfies
%\[
%\caL_0\left[\caL_0\tilde{\u}_2+H_0\partial_r\tilde{\u}_1-\frac{1}{2}\nabla_\u^3W(\phi_h)(\tilde{\u}_1,\tilde{\u}_1)+\nabla_\u{\bf V}(\phi_h)\tilde{\u}_1+rH_1\partial_r\phi_h\right]
%=-(\eta_1-\eta_2)\partial_r^2\phi_h+{\bf B}_2(\tau),
%\] 
%and 
%\[
%\mathbf{m}(\tau)\cdot \mathbf{M}=-M_1\frac{\int_S\left[ -|\nabla_s H_0|^2 +\frac{\eta_1+\eta_2}{2} H_0^2-\frac12 H_0^2(H_0-a_0)^2-H_1H_0(H_0-a_0) \right] \rmd s}{\int_S H_0^2\rmd s}.
%\]
 \end{Formal}

Finally, we break Assumption \ref{a:homoex}, and present a Birkhoff billiard potential for which the stable and unstable manifolds 
associated the origin of (\ref{e:leadode}) intersect non-transversely.  The result is a one-parameter family of homoclinics $\Phi_h(\cdot;\bm;\theta)$  which 
are distinct up to translation, see Example\,\ref{e:2} and the right-most image in  Figure\,\ref{f:billiard}. The corresponding linearized operators 
$\caL_{\delta,\theta}$ possess a two-dimensional kernel spanned by the translational eigenmode and $\partial_\theta\Phi_h(\cdot;\bm;\theta)$. The higher order 
kernel suggests a geometric evolution not only of the  shape of the underlying interface,  $\Gamma$, but also a possibility for a continuous evolution of the 
composition of the underlying bilayer profile, $\theta=\theta(s,t)$ where $s$ parameterizes position along the interface $\Gamma$. In such a system the 
geometric and compositional evolution could couple, leading to a model that couples geometric evolution with lipid composition, as depicted in the bottom-left 
panel of Figure\,\ref{f:packing}.

%%%%%%%%%%%%%%%%%%%%%%%%%%%%%%%%%%%%%%%%%%%%%%%%%%%%%%
\section{Quasi-bilayers: A reduction to the Canham-Helfrich free energy}\label{s:2}
\begin{figure}
  \begin{center}
    \includegraphics[width=0.5\textwidth]{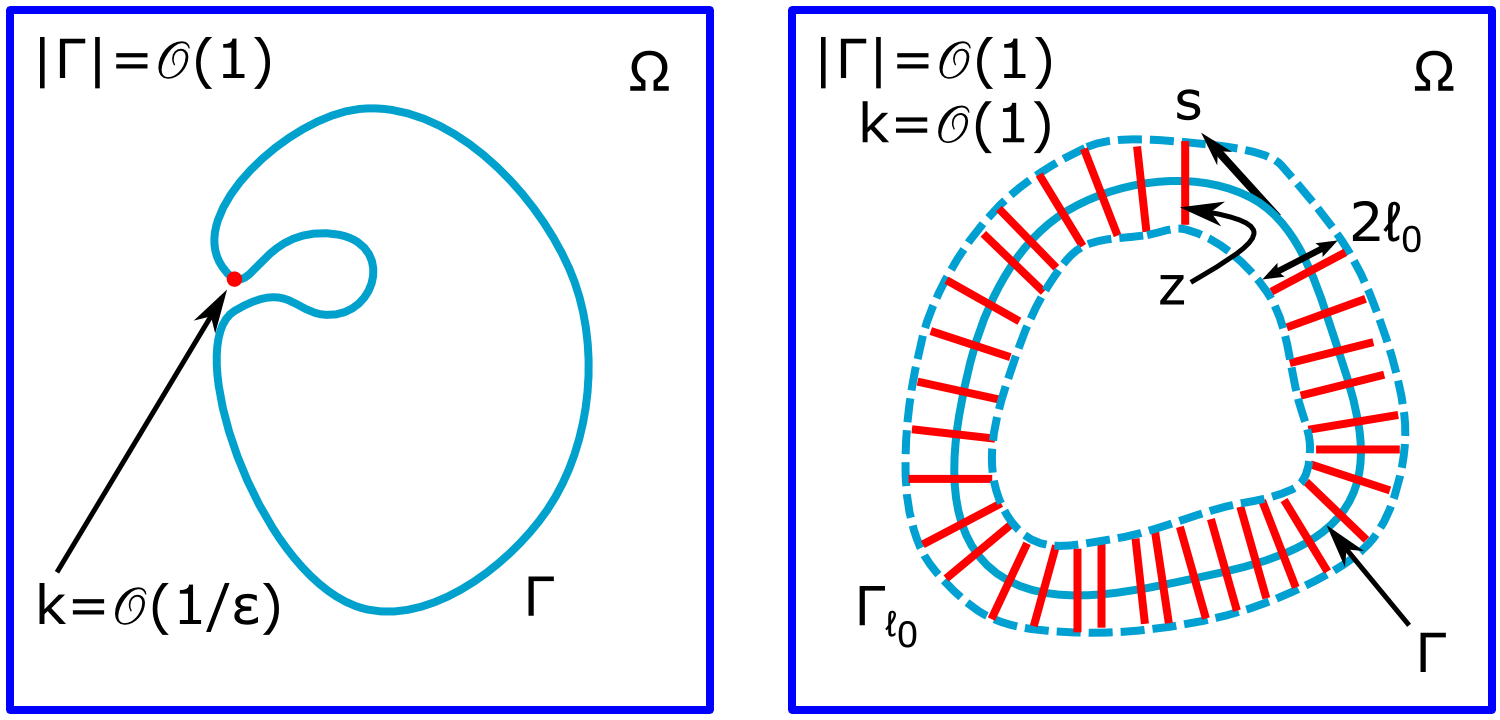}
  \end{center}
  \caption{In both panels, the blue rectangle represents the physical domain $\Omega$ of the amphiphilic mixture while the light blue curve represents the interface $\Gamma$ of the amphiphilic morphology. The left panel shows an inadmissible interface where local whisker coordinates are not well-defined, due to the existence of points (the red dot) where the curvature of the interface $k$ is of order $\caO(1/\veps)$. The right panel depicts an admissible interface $\Gamma$: Its length $|\Gamma|$ is of order $\caO(1)$ and its curvature is of order $\caO(1)$ everywhere. For such an interface, there exists a constant $l_0>0$ such that the union of $6l_0$-long whiskers (red lines) form a neighborhood of the interface, depicted as the annulus enclosed by light blue dash curves, denoted by $\Gamma_{l_0}$. In addition, no two distinct whiskers of length $6\ell_0$ intersect. For an admissible interface the whiskered coordinates 
 $(z,s)$ give a smooth change of variables of $\Gamma_{l_0}$, where $z$ is the $\veps$-scaled signed distance to the interface and $s$ is the tangential variable of the interface.}
\label{f:whisker}
\end{figure}

%%%%%%%%%%%%%%%%%%%%%%%%%%%%%%%%%%%%%%%%%%%%%%%%%%%%%%
The construction of quasi-bilayer profiles starts with a rigorous definition of an admissible interface $\Gamma$. Specifically, 
we assume the interface $\Gamma$ is a smooth $d-1$ dimensional manifold admitting a volume-preserving parameterization, 
$\Xi(s): S\subset \R^{d-1}\rightarrow \Gamma$. For every $s\in S$ and $l>0$, we can define a whisker $\omega(s,l)$ based at $\Xi(s)$, that is,
\[
\omega(s,l):=\{\Xi(s)+\veps z{\bf n}(s)\mid |z|\leq \frac{l}{\veps}\},
\]
where ${\bf n}(s)$ is the outward unit normal vector of the interface $\Gamma$ at point $\Xi(s)$ and $z$ is the $\veps$-scaled, signed distance to the interface $\Gamma$. We restricted ourselves to interfaces far from self-intersections.
\begin{Definition}[Admissible interfaces with reach $l_0$]\label{d:noselfin}
An admissible interface $\Gamma\subset\R^d$ is a smooth $d-1$ dimensional manifold far from self-intersections. More precisely, there exists $l_0>0$ such that 
\begin{itemize}
\item[(\rmnum{1})] No two whiskers of length $6\ell_0$ intersect, that is, $\omega(s_1,3l_0)\cap \omega(s_2,3l_0)=\emptyset$, for any $s_1\neq s_2\in S$.
\item[(\rmnum{2})] For each $l\in(0,3l_0]$ the set
\[
\Gamma_{l}:=\bigcup_{s\in S}\omega(s,l), 
\]
is a neighborhood of the interface $\Gamma$.
\end{itemize}
\end{Definition}
We call the set $\Gamma_{l_0}$ the \textit{inner region} and the set $\Omega\backslash \Gamma_{l_0}$ the \textit{outer region}, see Figure \ref{f:whisker}. In addition, we call $(z,s)\in[-\frac{l_0}{\veps},\frac{l_0}{\veps}]\times S\subset\R^d$ the \textit{local whisker coordinates}. The construction of the quasi-bilayer is based on a formal asymptotic analysis in the inner and outer regions; see Figure \ref{f:whisker}. Throughout this section, we fix $N=2$.

%%%%%%%%%%%%%%%
%\paragraph{Inner expansion}
%%%%%%%%%%%%%%%
To construct quasi-bilayer solutions of the full PDE \eqref{e:smfch} within the inner region $\Gamma_{l_0}$, we derive an ODE system whose homoclinic solutions have sufficient flexibility
to approximate the full behavior of the reduction of the PDE to the whiskers.  To begin, we transform the stationary mFCH equation \eqref{e:smfch} to the local whiskered coordinates $(z,s)$. 
In the whiskered coordinated the Laplacian becomes
\beq\label{e:lap}
\veps^2\Delta=\partial_z^2+\veps H(z,s)\partial_z+\veps^2\Delta_G,
\eeq
where at leading order the operator $\Delta_G$ reduces to the Laplace-Beltrami operator and $H$ is the extended curvature of the expression,
\beq
\label{e:H-def}
H(z,s)=-\sum_{j=1}^{d-1}\frac{k_j(s)}{1-\veps rk_j(s)}.
\eeq
%=H_0(s)+\caO(\veps), \text{where }H_0(s):=-\sum_{j=1}^{d-1}k_j(s),
Here $\{k_j(s)\}$ is the set of principal curvatures of $\Gamma$ at $\Xi(s)$; see  \cite{doelmanpromislow_2014} for details about $H$ and $\Delta_G$.
We plug the Laplacian \eqref{e:lap} into the stationary mFCH equation \eqref{e:smfch}, yielding, 
\begin{equation}
\label{e:innerpde}
\begin{aligned}
&\left(\partial_z^2+\veps H(z,s) \partial_z+\veps^2 \Delta_G-\nabla_\u^2W(\u)+\veps(\nabla_\u {\bf V}(\u))^{\rm T}\right)\left(\partial_z^2\u+\veps H(z,s)\partial_z\u+\veps^2\Delta_G\u-\nabla_\u W(\u)+\veps{\bf V}(\u)\right)\\
&=\veps^2\left[\eta_2\nabla_\u W(\u)-\eta_1(\partial_z^2 \u+\veps H(z,s)\partial_z\u+\veps^2\Delta_G\u)+{\bf m}\right],
\end{aligned}
\end{equation}
which, up to the leading order, reduces to the following ODE system
\[
(\partial_z^2-\nabla_\u^2 W(\u))(\partial_z^2\u-\nabla_\u W(\u))=0.
\]
As a starting point in our construction we make the following assumption.
\begin{Assumption}\label{a:homoex}
The leading-order Hamiltonian ODE system \eqref{e:leadode},
 \[
 \partial_z^2\u-\nabla_\u W(\u)=0,
 \]
  admits, up to translations, an orbit $\phi_h$ that is homoclinic to the origin. 
% We don't NEED this, we just LIKE it
%The profile $\phi_h$ is asymmetric in the sense that it satisfies the following symmetry-breaking condition,
%\beq\label{e:symbrek}
%\int_\R {\bf V}(\phi_h)\cdot \partial_z\phi_h \rmd z\neq 0.
%\eeq
Moreover, the linearized operator 
%\beq
%\label{e:L0}
%\begin{matrix}
%\caL_0: & (H^2(\R))^2 & \longrightarrow & (L^2(\R))^2\\
%& \u & \mapsto & \partial_z^2\u-\nabla_\u^2W(\phi_h)\u, 
%\end{matrix}
%\eeq
\beq
\label{e:L0}
\caL_0:=  \partial_z^2 -\nabla_\u^2W(\phi_h), 
\eeq
admits $0$ as a simple eigenvalue with the corresponding eigenspace spanned by $\phi_h^\prime(z)$.
\end{Assumption}
\begin{Remark}
A homoclinic orbit satisfying the condition 
\beq\label{e:symbrek}
\int_\R {\bf V}(\phi_h)\cdot \partial_z\phi_h \rmd z\neq 0.
\eeq 
is called asymmetric as it cannot have a translate that is even about $z=0$. We will discuss potentials $W$ which satisfy Assumption \ref{a:homoex} in Section \ref{s:3}.
\end{Remark}

To develop approximate solutions of  ({e:innerpde}) we assume $\u$ is independent of $s$ and replacing the extended curvature $H$ with a Melnikov 
parameter $a$, yielding the whiskered ODE system,
\begin{numcases}
{\label{e:2order}} \partial_z^2\u+\veps a\partial_z\u-\nabla_\u W(\u)+\veps{\bf V}(\u) = \veps^2 \v, \label{e:2order1}\\
(\partial_z^2+\veps a\partial_z-\nabla_\u^2W(\u)+\veps(\nabla_\u {\bf V}(\u))^{\rm T})\v=\eta_2\nabla_\u W(\u)-\eta_1(\partial_z^2\u+\veps a \partial_z\u)+{\bf m}.\label{e:2order2}
\end{numcases}
Assuming $\u(z)=\phi_h(z)+\caO(\veps)$ in \eqref{e:2order}, the leading order term of $\v$, denoted $\zeta_h(z)$, satisfies the equation, 
\beq
\caL_0\zeta_h(z)=(\eta_2-\eta_1)\nabla_\u W(\phi_h(z))+{\bf m}.
\eeq
It is natural to assume $\zeta_h\perp \ker(\caL_0)= \span\{\partial_z\phi_h\}$ since the projection of $\zeta_h$ onto the kernel of $\caL_0$ corresponds to a translation of $\zeta_h$ in $z$. With this assumption we have
\[
\zeta_{h}:=\caL_{0,\perp}^{-1}\left((\eta_2-\eta_1)\nabla_\u W(\phi_h)+{\bf m}\right),
\]
where $\caL_{0,\perp}^{-1}$ is the bounded inverse of the operator $\caL_0$ restricted to the orthogonal complement of  $\span\{ \partial_z\phi_h\}$ in $(L^2(\R))^2$. Noting that $\zeta_h$ converges to a constant state in infinity, that is,
\[
\lim_{z\rightarrow\pm\infty}\zeta_h(z)=-(\nabla_\u^2 W(0))^{-1}{\bf m}:={\bf E}
\]
Substituting ${\bf E}$ for $\v$ in \eqref{e:2order1},  yields the whisker ODE system
\beq
\label{e:melode}
 \partial_z^2\u+\veps a\partial_z\u-\nabla_\u W(\u)+\veps{\bf V}(\u) =\veps^2 {\bf E},
\eeq
which incorporates the Melnikov parameter $a$ and the far-field parameter ${\bf E}$.
We have the following lemma.
\begin{Lemma}\label{l:homoper}
Given Assumption \ref{a:homoex} and any ${\bf m}\in\R^2$, for sufficiently small $\veps>0$, there exists a unique choice of Melnikov parameter $a(\veps)$,
given by (\ref{e:aexp}) for which  
the whisker ODE system \eqref{e:melode},
%\[
% \partial_z^2\u+\veps a\partial_z\u-\nabla_\u W(\u)+\veps{\bf V}(\u) =\veps^2 {\bf E},
%\]
admits a homoclinic orbit 
\beq\label{e:homoper}
\Phi_{h}(z;\veps)=\phi_h(z)+\caO(\veps),
\eeq
connecting to the equilibrium 
\[
\Phi_{\infty}(\veps)=\veps^2{\bf B}+\caO(\veps^3),
\]
where the leading-order background state is given by $\mathbf{B}=(\nabla_\u^2W(0))^{-2}{\bf m}$. 
%In addition, the nonzero Melnikov parameter $a(\veps)$ admits the expansion as defined in \eqref{e:aexp}
%\[
% \quad a_0=-M_1^{-1}\int_\R  {\bf V}(\phi_h(z))\cdot \partial_z\phi_h(z) \rmd z, \quad M_1=\int_\R |\partial_z\phi_h(z)|^2\rmd z.
%\]
\end{Lemma}
\begin{Proof}
We rewrite the $4$th-order ODE system \eqref{e:melode} in first order form, that is,
\begin{equation}\label{e:1order}
\partial_z \U=F(\U,a,\veps),
\end{equation}
where $\U:=(\u,\v)$ and 
\[
F=\begin{pmatrix}
\v\\
\nabla_\u W(\u)-\veps{\bf V}(\u)-\veps a \v+\veps^2 {\bf E}
\end{pmatrix}.
\]
For $\veps=0$, the ODE system \eqref{e:1order} reduces to the leading order Hamiltonian ODE \eqref{e:leadode} in its first order form
\begin{equation}\label{e:0ode}
\partial_z\U=\begin{pmatrix}
\v\\
\nabla_\u W(\u)
\end{pmatrix},
\end{equation}
which, according to Assumption \ref{a:homoex}, admits a homoclinic orbit, $\U_h(z):=(\phi_h(z),\partial_z\phi_h(z))^{\rm T}$,
connecting to the hyperbolic equilibrium $0$. The persistence of the background equilibrium follows from the hyperbolicity of $W$ at the origin;
given any $a\in\R$ and sufficiently small $\veps>0$, there exists a hyperbolic equilibrium $\U_\infty(a,\veps)$ of the ODE \eqref{e:1order} with the expansion,
\[
\U_\infty(a, \veps)=\begin{pmatrix}\veps^2{\bf B}+\caO(\veps^3) \\ 0\end{pmatrix}.
\]
The persistence of the homoclinic solution $\U_h$ of the ODE \eqref{e:1order} follows via the Melnikov integral method. 
Linearizing \eqref{e:1order} at $\U=\U_h$ when $\veps=0$, yields the system
\beq\label{e:linho}
\partial_z \U(z)=\begin{pmatrix}
0 & {\bf I}_2 \\
\nabla_\u^2W(\phi_h(z)) & 0 
\end{pmatrix}\U(z),
\eeq
with the corresponding adjoint ODE system %of \eqref{e:linho} is
\beq\label{e:linad}
\partial_z \widetilde{\U}(z)=-\begin{pmatrix}
0 & \nabla_\u^2W(\phi_h(z)) \\
  {\bf I}_2 & 0 
\end{pmatrix}\widetilde{\U}(z),
\eeq
%due to the fact that $\nabla_\u^2W(\phi_h)$ is symmetric.
Denoting the stable and unstable manifold of the equilibrium $\U_\infty(a,\veps)$  as $\mathscr{M}^s(\U_\infty(a,\veps))$, $\mathscr{M}^u(\U_\infty(a,\veps))$ 
respectively, it is straightforward to see that
\beq
\begin{aligned}
T_{\U_h(z)}\mathscr{M}^s(0)\cap T_{\U_h(z)}\mathscr{M}^u(0)&=\span\{ \partial_z\U_h(z)\};\\
\left(T_{\U_h(z)}\mathscr{M}^s(0) \oplus T_{\U_h(z)}\mathscr{M}^u(0)\right)^\perp&=\span\{ \U_{ad}(z)\},
\end{aligned}
\eeq
where $\U_{ad}(z)=\left(-\partial_z^2\phi_h(z),\partial_z\phi_h(z)\right)^{\rm T}$.
Accordingly, we conclude that, for a small open ball $\mathscr{U}$ centered at $\U_h(0)$ in the hyperplane $\U_h(0)+\left(\span\{ \partial_z\U_h(0)\}\right)^\perp$, the two $2$-dimensional manifolds,
\[
\mathscr{U}\cap \left(\mathscr{M}^{s \backslash u}(0)\oplus \span\{ \U_{ad}(0)\}\right),
\]
 intersect transversally, forming a line, $\U_h(0)+\span\{ \U_{ad}(0) \}$. Moreover, the transversal intersection persists as we turn on the $\veps$-perturbation,
 that is, for sufficiently small $\veps$, the manifolds $\mathscr{U}\cap \left(\mathscr{M}^s (\U_\infty(a,\veps))\oplus \span\{ \U_{ad}(0)\}\right)$ and $\mathscr{U}\cap \left(\mathscr{M}^u(\U_\infty(a,\veps))\oplus \span\{ \U_{ad}(0)\}\right)$ also intersection transversally, and hence there exist two unique solutions $\U^s(z;a,\veps)$ defined on $z\in[0,\infty)$
and $\U^u(z;a,\veps)$ defined on $z\in(-\infty,0]$ of the full system \eqref{e:1order} such that
\[
\U^s(0;a,\veps)-\U^u(0;a,\veps)\in \span\{ \U_{ad}(0)\}.
\]
The homoclinic orbit $\U_h$ persists if $U^s(0;a,\veps)=\U^u(0;a,\veps)$ which in light of the relation above is equivalent to the Melnikov condition
\beq\label{e:meleq}
 \U_{ad}(0)\cdot\big[ \U^s(0;a,\veps)-\U^u(0;a,\veps)\big]=0.
\eeq
The left-hand side of the Melnikov condition admits the Taylor expansion
\begin{equation*}
\begin{aligned}
 \U_{ad}(0)\cdot\big[ \U^s(0;a,\veps)-\U^u(0;a,\veps)\big]&=\int_\R \U_{ad}(z)\cdot \partial_\veps F(\U_h(z), 0) \rmd z+\caO(\veps)\\
& = -\int_\R  \partial_z\phi_h(z) \cdot \big[ {\bf V}(\phi_h(z))+a\partial_z \phi_h(z)\big]\rmd z +\caO(\veps),
\end{aligned}
\end{equation*}
and solving for $a$ yields the expression (\ref{e:aexp}), and the persistence of the homoclinic $\U_h$ follows.
%Assuming the persistence of the homoclinic orbit,  we inner-product both sides of the first equation in \eqref{e:melode} with $\partial_z\phi_{h,\veps}$ in $(L^2(\R))^2$, yielding
%\begin{equation*}
%\begin{aligned}
%\veps a(\veps)\int_\R|\partial_z\phi_{h,\veps}|^2\rmd z=&\veps^2\int_\R\langle\psi_{h,\veps}, \partial_z\phi_{h,\veps}\rangle\rmd z-\int_\R\langle \partial_z^2\phi_{h,\veps}-\nabla_\phi W(\phi_{h,\veps})+\veps{\bf V}(\phi_{h,\veps}),\partial_z\phi_{h,\veps} \rangle\rmd z\\
%=&- |\partial_z\phi_{h,\veps}|^2\mid_{-\infty}^{+\infty}+W(\phi_{h,\veps})\mid_{-\infty}^{+\infty}+\veps\int_\R \langle -{\bf V}(\phi_{h,\veps})+\veps\psi_{h,\veps}, \partial_z\phi_{h,\veps} \rangle\rmd z\\
%=&\veps\int_\R \langle -{\bf V}(\phi_{h,\veps})+\veps\psi_{h,\veps}, \partial_z\phi_{h,\veps} \rangle\rmd z,
%\end{aligned}
%\end{equation*}
%which naturally leads to the expression of $a(\veps)$ as in \eqref{e:aveps}. The fact that $a_0\neq 0$ is due to the asymmetric condition \eqref{e:symbrek} in Assumption \ref{a:homoex}. The proof the lemma now boils down to show the persistence of the homoclinic orbit $(\phi_h,\zeta_h)$. \textbf{TBC}.
\end{Proof}

Within the inner region $\Gamma_{l_0}$, the quasi-bilayer profile $\u_q$ equals $\Phi_h(z;\veps)$, that is,
\beq\label{e:quabi1}
\u_q(x;\veps):=\Phi_{h}(z(x);\veps)=\phi_h(z(x))+\veps\phi_{h,1}(z(x))+\veps^2\phi_{h,2}(z(x))+\caO(\veps^3) ,
\eeq
where $x\in\Gamma_{l_0}$ and $z\in[-l_0/\veps,l_0/\veps]$ and the error terms are comprised of smooth perturbations which are $\caO(\veps^3)$ in the $L^\infty$ norm. 
We plug the expansion \eqref{e:quabi1} into \eqref{e:melode} and evaluate its terms of order $\veps$ and $\veps^2$, yielding respectively
\[
\begin{cases}
\phi_{h,1}=&-\caL_{0,\perp}^{-1}({\bf V}(\phi_h)+a_0\partial_z\phi_h),\\
\phi_{h,2}=&\caL_{0,\perp}^{-1}\left(\zeta_h-a^\prime(0)\partial_z\phi_h-a_0\partial_z\phi_{h,1}+\frac{1}{2}\nabla_\u^3 W(\phi_h)(\phi_{h,1},\phi_{h,1})-\nabla_u{\bf V}(\phi_h)\phi_{h,1}\right).
\end{cases}
\]
From Lemma \ref{l:homoper}, in the inner region $\Gamma_{l_0}$, there exists a constant $c>0$ independent of $\veps$ so that
\begin{equation*}
\begin{aligned}
\veps^2\Delta \mathbf{u}_q-\nabla_\u W(\u_q)+\veps{\bf V}(\u_q)=&(\partial_z^2+\veps H(z,s)\partial_z)\Phi_h-\nabla_\u W(\Phi_h)+\veps{\bf V}(\Phi_h)\\
=&\veps(H(z,s)-a(\veps))\partial_z\Phi_h+\veps^2{\bf E}+\caO(\veps^3)\\
=&\veps(H_0(s)-a_0)\partial_z\phi_h+\veps^2{\bf E}+\caO(\veps^2\rme^{-c|z|}+\veps^3),\\
\end{aligned}
\end{equation*}
and 
\begin{equation*}
\begin{aligned}
\frac{1}{2}\veps^2\eta_1|\nabla \mathbf{u}_q|^2+\eta_2W(\bf{u}_q)=&\frac{1}{2}\eta_1|\partial_z\Phi_h|^2+\eta_2W(\Phi_h)
=\frac{1}{2}(\eta_1+\eta_2)|\partial_z\phi_{h}|^2 +\caO(\veps\rme^{-c|z|}+\veps^2),
\end{aligned}
\end{equation*}
where $H(r,s)=H_0(s)+\veps rH_1(s)+\caO(\veps^2)$, all the error estimates are in the $L^\infty$ norm and in the second estimate, we apply the fact that $\frac{1}{2}|\partial_z\phi_h|^2=W(\phi_h)$ for all $r\in\R$.
The Jacobian $J$ of the change of coordinates from $x$ into $(z,s)$ admits the expansion 
\beq\label{e:jacobian}
J(z,s)=\veps+\veps^2zH_0(s)+\caO(\veps^3);
\eeq
(see the appendix of \cite{hayprom_2014}). Applying the two estimates above, 
we evaluate the mFCH free energy evaluated at $\u_q$ restricted to $\Gamma_{l_0}$, obtaining
\beq
\begin{aligned}
\mathcal{F}_{M,inner}(\u_q):=&\int_{\Gamma_{l_0}} \frac{1}{2}|\veps^2\Delta \mathbf{u}_q-\nabla_\u W(\u_q)+\veps{\bf V}(\u_q)|^2-\veps^2\left(\frac{1}{2}\veps^2\eta_1|\nabla \mathbf{u}_q|^2+\eta_2W(\bf{u}_q)\right)\rmd x\\
=&\veps^3\int_S\int_{-l_0/\veps}^{l_0/\veps}\left[\frac{1}{2}|(H_0(s)-a_0)\partial_z\phi_h+\veps{\bf E}|^2-\frac{1}{2}(\eta_1+\eta_2)|\partial_z\phi_h(z)|^2\right]\rmd z \rmd s+\caO(\veps^4)\\
=& \frac{1}{2}\veps^3M_1\int_S\left[(H_0(s)-a_0)^2-(\eta_1+\eta_2)\right]\rmd s+\caO(\veps^4).
\end{aligned}
\eeq
In the last step, we use the fact that
\[
\int_\R \partial_z\phi_h(z) \cdot {\bf E} \rmd z=0, \quad \int_{-l_0/\veps}^{l_0/\veps}\rme^{-cz}\rmd z=\caO(1), \quad \int_{-l_0/\veps}^{l_0/\veps}1\rmd z=\caO(\veps^{-1}).
\]
The total mass vector of amphiphilies within the inner region $\Gamma_{l_0}$, denoted $\mathcal{M}_{inner}$, admits the expansion
\beq
\begin{aligned}
\mathcal{M}_{inner}(\u_q):=&\int_{\Gamma_{l_0}}\u_q\rmd x=\veps\int_S\int_{-l_0/\veps}^{l_0/\veps}\Phi_h(z;\veps)\left(1+\veps z H_0(s)+\caO(\veps^2)\right)\rmd z\rmd s\\
=&\veps\int_S1\rmd s\int_{-l_0/\veps}^{l_0/\veps}(\phi_{h}(z)+\veps\phi_{h,1}(z)+\veps^2{\bf B})\rmd z+\veps^2\int_SH_0(s)\rmd s\int_\R z\phi_h(z)\rmd z+\caO(\veps^3)\\
=&\veps|\Gamma|\int_\R\phi_{h}(z)\rmd z+\veps^2\left[|\Gamma|\left(2l_0{\bf B}+\int_\R\phi_{h,1}(z)\rmd z\right)+\int_SH_0(s)\rmd s\int_\R z\phi_h(z)\rmd z\right]+\caO(\veps^3),\\
\end{aligned}
\eeq
where we use the fact that the area of $\Gamma$ and $S$ are the same, i.e., $|\Gamma|=|S|$.

%%%%%%%%%%%%%%%
%\paragraph{Outer expansion}
%%%%%%%%%%%%%%%
Within the outer region $\Omega\backslash \Gamma_{l_0}$, the quasi-bilayer profile $\u_q$ makes a smooth transition to the asymptotic state of $\Phi_\infty(\veps,{\bf m})$ 
as specified in (\ref{e:quasib}). 
%To make the quasi-bilayer profile smooth, we need a smooth transition between the inner and outer region. To be precise, given a smooth cut-off function$\chi:\R\to\R$ so that $\chi(r)=1$ for $|r|\leq 1$ and $\chi(r)=0$ for $|r|\geq 3$, we define
%\begin{equation}
%\u_q(x;\veps):=
%\begin{cases}
%\Phi_\infty(\veps), & x\in\Omega\backslash \Gamma_{3l_0},\\
%(1-\chi(\frac{\veps |z(x)|}{l_0}))\Phi_\infty(\veps)+\chi(\frac{\veps |z(x)|}{l_0})\Phi_{h,\delta}(z(x)), & x\in \Gamma_{3l_0}\backslash \Gamma_{l_0},
%\end{cases}
%\end{equation}
In particular for any $x\in\Omega\backslash \Gamma_{l_0}$ we have the estimate,
\beq\label{e:interm}
\u_q(x)=\Phi_\infty(\veps)+\caO(\rme^{-c/\veps})=\veps^2{\bf B}+\caO(\veps^3).
\eeq
%Again, as in the inner region, the $ \caO(\veps^3)$ term represents a smooth perturbation of order $\veps^3$.
With these estimates we find that the outer region contribution to the mFCH free energy and the total mass vector of $\u_q$ are, respectively, 
%restricted to the outer region respectively as $\mathcal{F}_{M,outer}(\u_q)$ and $\mathcal{M}_{outer}(\u_q)$, we readily obtain
\beq
\begin{aligned}
\mathcal{F}_{M,outer}(\u_q):=&\int_{\Omega\backslash\Gamma_{l_0}} \frac{1}{2}|\veps^2\Delta \mathbf{u}_q-\nabla_\u W(\u_q)+\veps{\bf V}(\u_q)|^2-\veps^2\left(\frac{1}{2}\veps^2\eta_1|\nabla \mathbf{u}_q|^2+\eta_2W(\bf{u}_q)\right)\rmd x\\
=&(|\Omega|-2l_0|\Gamma|)\left[\frac{1}{2}|-\nabla_\u W(\Phi_\infty)+\veps{\bf V}(\Phi_\infty)|^2-\veps^2\eta_2W(\Phi_\infty)\right]+\caO(\rme^{-c/\veps})\\
=&\caO(\veps^4),\\
\mathcal{M}_{outer}(\u_q):=&\int_{\Omega\backslash\Gamma_{l_0}}\u_q\rmd x=\veps^2(|\Omega|-2l_0|\Gamma|){\bf B}+\caO(\veps^3).
\end{aligned}
\eeq

%Combining the inner and outer contributions, the quasi-bilayer profile given in (\ref{e:quasi}) satisfies 
%\begin{equation*}
%\u_q(x;\veps)=\begin{cases}
%\Phi_h(z(x),\veps), & x\in\Gamma_{l_0},\\
%(1-\chi(\frac{\veps |z(x)|}{l_0}))\Phi_\infty(\veps)+\chi(\frac{\veps |z(x)|}{l_0})\Phi_{h,\delta}(z(x)), & x\in \Gamma_{3l_0}\backslash \Gamma_{l_0},\\
%\Phi_\infty(\veps), & x\in\Omega\backslash \Gamma_{3l_0},
%\end{cases}
%\end{equation*}
has mFCH free energy $\mathcal{F}_M(\u_q)=\mathcal{F}_{M,inner}(\u_q)+\mathcal{F}_{M,outer}(\u_q)$ and $\veps$-scaled total mass vector 
$\mathcal{M}(\u_q)=\frac{1}{\veps}\left(\mathcal{M}_{inner}(\u_q)+\mathcal{M}_{outer}(\u_q)\right)$ which admit the following expansions
\begin{equation*}
\begin{cases}
\mathcal{F}_M(\u_q)=\veps^3\frac{M_1}{2}\int_S\left[(H_0(s)-a_0)^2-(\eta_1+\eta_2)\right]\rmd s+\caO(\veps^4),\\
\mathcal{M}(\u_q)= |\Gamma|\int_\R\phi_{h}(z)\rmd z+\veps\left[|\Omega|{\bf B}+|\Gamma|\int_\R\phi_{h,1}(z)\rmd z+\int_SH_0(s)\rmd s\int_\R z\phi_h(z)\rmd z\right]+\caO(\veps^2).
\end{cases}
\end{equation*}
These results establish Theorem \ref{t:canhel}.

%%%%%%%%%%%%%%%%%%%%%%%%%%%%%%%%%%%%%%%%%%%%%%%%%%%%%%
\section{Regularized-billiard potentials: Existence of homoclinics}\label{s:3}
%%%%%%%%%%%%%%%%%%%%%%%%%%%%%%%%%%%%%%%%%%%%%%%%%%%%%%
In this section, we introduce the Birkhoff-billiard potentials, $B:\R^N\mapsto \R$, and their regularized forms, denoted by $\{B_\delta\}$, where $\delta>0$ is a 
small free parameter independent of $\veps$.  We construction $n$-collision homoclinic orbits, $\psi$, of the Birkhoff-Hamiltonian ODE
\beq
\label{e:billode}
\partial_z^2\u-\nabla_\u B(\u)=0,
\eeq
and demonstrate their persistence under the regularization, yielding a homoclinic solution $\psi_\delta$ of
\beq
\label{e:quasiode}
\partial_z^2\u-\nabla_\u B_\delta(\u)=0,
\eeq
that satisfies Assumption \ref{a:homoex}. We perform a spectral analysis of the associated linearized operator 
\beq
\label{e:lin-op}
\caL_\delta=\partial_z^2-\nabla_\u^2B_\delta(\psi_\delta)
\eeq
showing that each collision in corresponds to a collision eigenvalue --  a positive, order of $\delta^{-2}$ eigenvalue of $\caL_\delta$. 
While the essential spectrum of $\caL_\delta$ is uniformly in the left-half complex plane, we give an example of a problem with 
hidden symmetries that generate a non-trivial kernel, in addition to the translational symmetry of the orbit $\psi_\delta$,
which leads to novel center-stable modes other than the meandering mode. Throughout this section, we restrict our attention to
the two-species case, $N=2$, for simplicity.
% and $\Gamma$ is flat, where the Laplacian \eqref{e:lap} becomes $\veps^2\Delta=\partial_{rr}+\veps^2\partial_{ss}$.  

%%%%%%%%%%%%%%%%%%%%%%%%%%%%%%%%%%%%
\subsection{Singular Hamiltonian: Homoclinics in a billiard limit potential}
%%%%%%%%%%%%%%%%%%%%%%%%%%%%%%%%%%%%
 The Birkhoff-billiard potential is defined as follows.
\begin{Definition}\label{d:bill}
A map $B:\R^2\mapsto \R$ defined on the physical domain,
\[
\Sigma:=\{{\bf u}\in\R^{2}\mid  u_1+u_2\leq 1, u_1, u_2\geq 0\},
\]
is a \textbf{Birkhoff billiard potential}, if it satisfies the following conditions.
\begin{itemize}
\item[(\rmnum{1})]There exists $R_0\in(0,1)$ such that the potential $B(\u)$ is radially symmetric in 
$\u$ on the quadrant $\mathscr{Q}:=\{\u\in\Sigma\mid |\u|<R_0\}$. In addition,
there exists an integer $ \ell\geq 2$, a positive constant $b^+$ and a strictly increasing $C^{\ell+1}$-smooth function 
$b(r)$ defined on $[0, R_0]$ with 
\begin{equation*}
b(r)=\begin{cases}
r^2, & r\in[0,R_0/3],\\
b^+, & r=R_0,
\end{cases}
\end{equation*}
so that  $B(\u)=b(|\u|)$, for $\u\in\mathscr{Q}$.
\item[(\rmnum{2})] The potential $B(\u)$ is piecewise constant on the \textbf{billiard region} $\scb:=\Sigma\backslash\mathscr{Q}$.
More specifically, there exist two mutually disjoint simply connected open sets $\scb^\pm\subset \scb$ with $\scb=\overline{\scb^+}\cup\overline{\scb^-}$ such that
\begin{equation}
B(\u)=
\begin{cases}
b_+, &\u\in\overline{\scb^+},\\
-b_-, &\u\in\scb^-.\\
\end{cases}
\end{equation}
In addition, the potential $B(\u)$ is $C^{\ell+1}$-smooth in the simply connected domain $\Sigma\backslash\scb^-$ and the common boundary, called the {\em collision curve}, $\sca:=\overline{\scb^+}\cap\overline{\scb^-}$ is $C^{\ell+1}$-smooth with positive length. There also exists a neighborhood of $\sca$ in the physical domain $\Sigma$, denoted as $\scn(\sca)$, and a $C^{\ell+1}$-smooth signed distance function 
$\rho(\u)$ defined on $\scn(\sca)$ such that 
\[
\rho(\u)\mid_{\u\in\scb^+\cap\scn(\sca)}>0, \hspace{0.6cm} \nabla_\u\rho(\u)\mid_{\u\in\scn(\sca)}\neq 0,   
\]
and the collision curve $\sca$ is the zero level set of $\rho(\u)$, that is,
\[
\sca\subseteq\{\u\in\scn(\sca)\mid \rho(\u)=0\}.
\]
\end{itemize}
\end{Definition}

We study the leading-order Hamiltonian ODE system \eqref{e:billode} 
%when the potential is a billiard potential, that is, 
%\beq
%\label{e:billode}
%\partial_z^2\u-\nabla_\u B(\u)=0,
%\eeq
%For conveniences, 
which we rewrite as a first-order system,
\beq
\label{e:ham}
\partial_z\begin{pmatrix}\u\\ \v\end{pmatrix}=\begin{pmatrix}\v\\ \nabla_\u B(\u)\end{pmatrix},
\eeq
with the Hamiltonian $H(\U)=\frac{1}{2}|\v|^2-B(\u)$ where $\U:=(\u,\v)^{\rm T}$.
For clarity, we use the term \textit{orbit} for an orbit in the phase space and the term \textit{trajectory} for the projection of an orbit onto the $\u$-plane. A point $\U=(\u,\v)^{\rm T}$ in the phase space is called an \textit{inner point} if $\u\in\scb^+$ and a \textit{collision point} if $\u\in \sca\cap \scb^o$. At the moment of a collision, say $z$, the flow of \eqref{e:ham} satisfies the reflection law,
that is, the angle of incidence equals the angle of reflection, or more precisely,
\beq\label{e:reflection}
{\bf n}\cdot \big(\v_+(z)+\v_-(z)\big)=0, \quad \v_+(z)-\v_-(z)=\left[{\bf n}\cdot \big(\v_+(z)+\v_-(z)\big) \right]{\bf n},
\eeq
where $\displaystyle\v_\pm(z):=\lim_{s\rightarrow z^\pm}\v(s)$ and ${\bf n}$ is the unit normal vector of the collision curve $\sca$ at point $\u(z)$. 
\begin{Remark}
The Birkhoff-billiard potentials have been studied when the collision curve $\sca$ is piecewise smooth, 
however the reflection law fails when the normal vector is not well-defined.  Such a collision is generically not well-defined 
unless the billiard region $\scb^+$ is a fundamental domain of a finite Coxeter group, see \cite{treshchev_1991} for details. 
\end{Remark}

%For conveniences, the term ``orbit" refers to any orbit in the phase space and the term ``trajectory" refers to the projection of an orbit to the physical configuration space- the $\u$-plane. 
Our assumption that the origin $\u=0$ is a strict local minimum point of $B(\u)$ indicates that the origin is an equilibrium for the Hamiltonian ODE \eqref{e:ham}, 
admitting a two dimensional stable manifold,  $\mathscr{M}^s(0)$, and a two dimensional unstable manifold, $\mathscr{M}^u(0)$.  
While $\mathscr{M}^s(0)$ and  $\mathscr{M}^u(0)$ lie in the three dimensional invariant manifold $\mathscr{H}_0:=\{\U\in\R^4\mid H(\U)=0\}$, the conditions of 
Definition \ref{d:bill} are not sufficient to conclude that the intersection $\mathscr{M}^s(0)\cap \mathscr{M}^u(0)$ gives an orbit homoclinic to the origin. 
We introduce two main assumptions--the existence and transversality of a homoclinic orbit.
\begin{Assumption}[\textbf{Existence of a homoclinic orbit}]\label{a:homo}
There exists $n\in\Z^+$ and $z_1<z_2<\cdots<z_n$ such that 
the Hamiltonian system \eqref{e:ham},
\[
\partial_z\begin{pmatrix}\u\\ \v\end{pmatrix}=\begin{pmatrix}\v\\ \nabla_\u B(\u)\end{pmatrix},
\]
 with a billiard potential $B(\u)$
 admits an $n$-collision orbit homoclinic to the origin, denoted by $\Psi_h(z)=\big(\psi_h(z),\partial_z\psi_h(z)\big)^{\rm T}$, satisfying that
\begin{subnumcases}{\label{h:main2} }
\text{\rm $n$-collision}: \{z\in\R\mid \psi_h(z)\in \sca\}=\{z_i\}_{i=1}^n\subset \scb^o, &\label{h:1}\\
\text{\rm piecewise-linearity}:  \{\psi_h(z)\mid z\in(z_i, z_{i+1})\}\subset \scb^+,  i=1,2,\cdots, n-1,&\label{h:3}\\
\text{\rm Non-tangency}: {\bf n}_i\cdot\partial_z\psi_h(z_i)\neq 0,  i=1,2,\cdots, n, &\label{h:2}
\end{subnumcases}
 where ${\bf n}_i$ is the unit normal vector of $\sca$ at the $i$-th collision point ${\bf c}_i:=\psi_h(z_i)$ for all $i$. Without loss of generality, we set $z_1=0$.
 % In addition, the linearized operator
% \beq
% \label{e:L}
% \begin{matrix}
% \caL_h: & (H^2(\R))^2 & \longrightarrow & (L^2(\R))^2\\
% &\v(z)&\longmapsto&\v_{rr}-\nabla_\u_2W(\u_h)\v, 
% \end{matrix}
% \eeq
% admits $0$ as its simple eigenvalue and $\span\{\dot{\u}_h\}$ as the corresponding eigenspace.
\end{Assumption}
\begin{figure}
\centering
\includegraphics[scale=0.28]{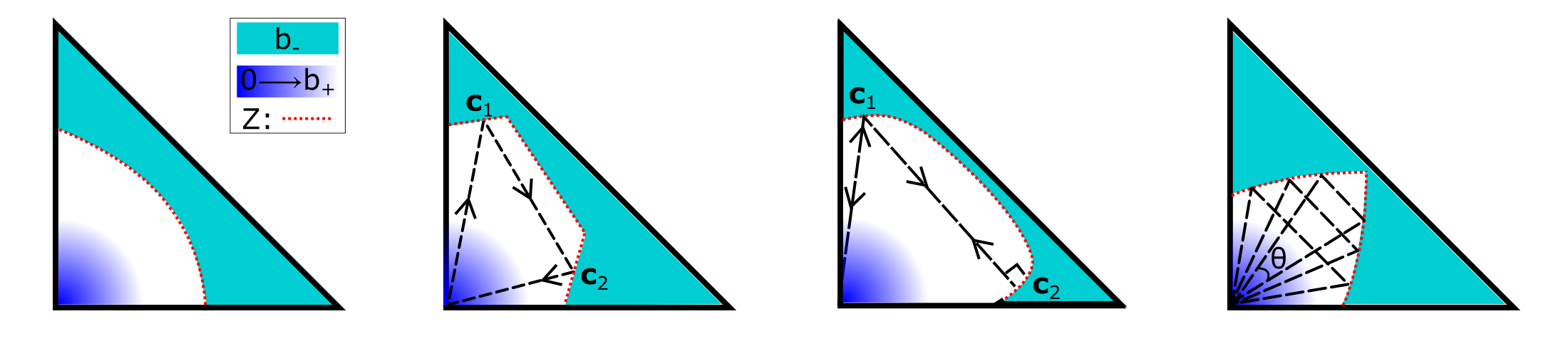}
\caption{In each panel, the triangle denotes the physical domain $\Sigma$, the potential is positive by smaller than $b^+$ in the faded blue region, the white region is $\scb^+$,
the cyan region is $\scb^-$,  the red dotted line in-between is the collision curve, $\sca$, while the black dashed lines denote trajectories of a homoclinic orbit. 
The leftmost panel gives a billiard potential. The middle two panels present a 2-collision homoclinic (middle left) and a 3-collision homoclinic (middle right), both admitting two 
distinct collision points. The middle left panel presents a billiard potential satisfying both Assumption \ref{a:homo} and  \ref{a:trans}. The rightmost panel is a {\em universal billiard} potential 
which admits a family of homoclinics parameterized by the angle $\theta$. This potential shows that symmetries other than translation may exist; 
see Example \ref{e:2} for more details.}
\label{f:billiard}
\end{figure}

The condition (\ref{h:1}) implies that the homoclinic trajectory intersect  the collision curve $\sca$ precisely $n$ distinct moments, but
not necessarily $n$ distinct points; see Figure \ref{f:billiard}. The condition (\ref{h:3}) requires the homoclinic orbit does not leave $\scb_+$ between collisions,
and hence the orbit is piecewise linear between collisions. The condition \eqref{h:2} guarantees that the collisions are \textit{not tangential}: a tangential collision, ${\bf n}_i\cdot\partial_z\psi_h(z_i)=0$, leads to a degeneracy--a loss of the smoothness of the flow and a generic failure of the persistence of the homoclinic orbit with respect to perturbations, see \cite{turaev_2007} for details about degenerate collisions. Figure \ref{f:billiard} presents two generic cases with direct physical interpretations: a 2-collision homoclinic orbit (middle left panel) 
corresponds to an asymmetric $AB$-type bilayer, and a 3-collision homoclinic with 2 distinct collision points (middle right panel) corresponding
to a symmetric $ABA$-type bilayer, where $A$ and $B$ are the phase of two distinctive amphiphilies, see also Figure \ref{f:SBlayer}. 

 Assumption \ref{a:homo} can be readily satisfied by a straightforward tuning of the collision curve $\sca$, but not all homoclinics 
 persist under smooth perturbations of the billiard potential. The following is an important class of counterexample that admits a one parameter
 family of co-existing homoclinics.
\begin{Example}[{\bf A universal Birkhoff-billiard potential}]\label{e:2}
For any fixed $c\in(0,\frac{7}{8}]$, we define the universal billiard potential,
\begin{equation*}
\Theta_c(\u)=\begin{cases}
|\u|^2, & |\u|\leq 1/4,\\
\frac{1}{4}, &\u\in\{\u\in\Sigma\mid  |\u|+\frac{1}{\sqrt{2}}|u_1-u_2|\leq c, |\u|>\frac{1}{2} \}, \\ 
-\frac{1}{4}, & \u\in\{\u\in\Sigma\mid  |\u|+\frac{1}{\sqrt{2}}|u_1-u_2|> c\},\\
\end{cases}
\end{equation*}
yielding, up to translation, a family of homoclinic orbits connecting the origin, parameterized by the angle $\theta$ of the homoclinic trajectory formed at the origin; see the rightmost panel of Figure \ref{f:billiard}. We anticipate the existence of a smoothing of the universal billiard potential, denoted as $\Theta_{c,\delta}$ which preserves a large part of the family of homoclinics. 
As a result, the dimension of the kernel of the linearized operator $\caL_\delta:=\partial_z^2-\nabla_\u^2\Theta_{c,\delta}$ will be larger than one.
\end{Example}

In order to guarantee the persistence of homoclinic orbits under smooth perturbations of the potential, 
we impose the following transversality assumption.
\begin{Assumption}[\textbf{Transversality of the homoclinic orbit}]\label{a:trans}
The homoclinic orbit $\Psi_h$ of the ODE \eqref{e:ham} comprises the transversal intersection of the stable and unstable manifold of the origin, that is, for all $z\in\R$,
\[
{\bf T}_{\Psi_h(z)}\mathscr{M}^s(0)\cap {\bf T}_{\Psi_h(z)}\mathscr{M}^u(0)=\span\{\Psi_h^\prime(z)\}, 
\]
where ${\bf T}_\mathbf{p}\mathscr{M}$ represents the tangential space of the point $\mathbf{p}$ in the manifold $\mathscr{M}$.
\end{Assumption}

\begin{Remark}
The set of billiard potentials satisfying Assumption \ref{a:homo} and \ref{a:trans} is not empty. In fact, families of such potentials can be readily 
constructed by straightforward ray-tracing: Given arbitrary two distinctive points ${\bf c}_1,{\bf c}_2\in\Sigma\backslash\{0\}$,  we can tune 
the boundary $\sca$ according to the reflection law \eqref{e:reflection} so that the triangle  with its vertex ${\bf c}_1,{\bf c}_2$ and the origin 
is a unique trajectory of the ODE $\eqref{e:ham}$ that is homoclinic to the origin; see the middle left panel of Figure \ref{f:billiard}.
\end{Remark}

%\begin{Remark}
%The Assumption \ref{a:trans} is generic but not universal. It is quite reasonable to construct a billiard potential with hidden symmetries which admits a continuum of homoclinics and thus the multiplicity of eigenvalue $0$ of the corresponding $\caL_\delta$ is big than one. The universal billiard constructed in Example \ref{e:2} is a concrete example where the multiplicity of eigenvalue $0$ is at least 2 and the 2 collisions still gives rise to 2 collision eigenvalues.
%\end{Remark}

%\begin{Example}[{\bf No homoclinic orbit in $\Sigma$}]\label{e:1}
%%An example where $\mathscr{M}^s(0)\cap \mathscr{M}^u(0)$ is only the origin.
%We choose a billiard potential satisfying
%\begin{equation*}
%W(\u)=\begin{cases}
%|\u|^2, & |\u|\leq 1/4,\\
%\frac{1}{4}, &\u\in\{\u\in\Sigma\mid u_2\leq 3/4\}\backslash\{\u\in\Sigma\mid |\u|\leq 1/2\}, \\ 
%-\frac{1}{4}, & \u\in\{\u\in\Sigma\mid u_2> 3/4\}.\\
%\end{cases}
%\end{equation*}
%It is straightforward to conclude that for this specific choice of potential, there is no homoclinic orbit in $\Sigma$ connecting to the origin.
%\end{Example}

%%%%%%%%%%%%%%%%%%%%%%%%%%%%%%%%%%%%
\subsection{Nonsingular Hamiltonian: Homoclinics in a regularized-billiard potential}
%%%%%%%%%%%%%%%%%%%%%%%%%%%%%%%%%%%%
The homoclinic orbit in the Birkhoff-billiard potential is singular in the sense that
 its first derivative admits a jump at each collision due to the reflection law. Being continuous but not $C^\ell$-smooth 
for any $\ell>0$ makes linearization difficult. Consequently we introduce the
 \textit{regularized Birkhoff-billiard potential}, whose dynamics retains the major features of non-smooth Birkhoff-billiard system.
In this sub-section we establish the persistence of transverse homoclinic orbits under the regularization. 
The dynamics away from the collision curve $\sca$ are smooth and easy to manipulate; the delicate part is to
 smoothen the billiard potential near the boundary $\sca$ while qualitatively minimizing the change of the dynamics nearby. 
%the set of corner points $\Gamma$, denoted as $\scn(\Gamma)$,  so that $\u_1, \u_2\notin \scn(\Gamma)$. In addition, we note that, according to the definition of $\Gamma$, any arc of the boundary $Z$ that connects two neighboring corner points is $C^{\ell+1}$-smooth. We name these arcs the \textit{boundary arcs} and label them as $\{\Gamma_i\}_{i\in I}$. For each boundary arc $\Z_i$, we fix a sufficiently small neighborhood in $\scb$, denoted as $\scn(\Gamma_i)$,  so that
%\[
%\scn_i\cap\scn_j=\emptyset, \quad i\neq j,
%\]
%where $\scn_i:=\scn(\Gamma_i)\backslash\scn(\Gamma)$. 

To introduce regularized-billiard potentials, we define a smooth jump function $\widetilde{\chi}_{\delta}(\rho):=(\widetilde{\chi}\ast h_\delta)(\rho)$,
where
\[
\widetilde{\chi}(\rho)=
\begin{cases}
b_+, & \rho\geq0,\\
-b_-,&\rho<0.
\end{cases}
\]
and $h$ is the classical mollifier with $h_\delta(\rho):=\delta^{-1}h(\delta^{-1}\rho)$ for $\delta>0$. 
%Here the mollifier $h$ takes the form
%\[
%h(\rho)=\begin{cases}
%c\exp(\frac{1}{\rho^2-1}), & |\rho|<1,\\
%0, & |\rho|\geq 1,
%\end{cases}
%\]
%where $c:=\int_\R\exp(\frac{1}{\rho^2-1})\rmd \rho$.
There exists $\delta_0>0$ so that $[-\delta_0,\delta_0]$ is part of the range of $\rho$. For any $\delta\in[-\delta_0,\delta_0]$, we introduce the \textit{smooth region}
\[
\scn_\delta(\sca):=\{\u\in\scn(\sca)\mid-\delta\leq\rho(\u)\leq\delta\}.
\]
and the regularized Birkhoff-billiard potential and refer to \cite{turaev_1998, turaev_2007} for a more general definition.
\begin{Definition}\label{d:quasi}
Given a billiard potential $B(\u)$ and any $\delta\in(0,\delta_0]$, the potential
\beq
B_\delta(\u)=
\begin{cases}
\widetilde{\chi}_{\delta}(\rho(\u)), & \u\in\scn_\delta(\sca),\\
B(\u),&\u\notin\scn_\delta(\sca).
\end{cases}
\eeq 
is a regularized-billiard potential with respect to the billiard potential $B(\u)$.
\end{Definition}

Replacing the Birkhoff-billiard potential with the regularization yields a smooth Hamiltonian system
\beq
\label{e:qham}
\begin{pmatrix}\dot{\u}\\\dot{\v}\end{pmatrix}=\begin{pmatrix}\v\\ \nabla_\u B_\delta(\u)\end{pmatrix},
\eeq
with Hamiltonian $H(\U; \delta)=\frac{1}{2}|\v|^2-B_\delta(\u)$. By definition \ref{d:quasi}, the origin is again a hyperbolic equilibrium of the system \eqref{e:qham} admitting
a two-dimensional stable manifold, $\mathscr{M}^s(0;\delta)$, and a two-dimensional unstable manifold, $\mathscr{M}^u(0;\delta)$.
We claim that, for sufficiently small $\delta$, transverse homoclinic orbits in the limit billiard system \eqref{e:ham} persist as solutions of \eqref{e:qham}. 
\begin{Proposition}\label{p:per}
Given a Birkhoff billiard potential $B(\u)$ for which the ODE \eqref{e:ham},
%\[
%\partial_z\begin{pmatrix}\u\\ \v\end{pmatrix}=\begin{pmatrix}\v\\ \nabla_\u B(\u)\end{pmatrix},
%\]
 admits an $n$-collision homoclinic $\Psi_h=(\psi_h,\partial_z\psi_h)^{\rm T}$ satisfying Assumptions \ref{a:homo} and \ref{a:trans}, then
for any sufficiently small $\delta>0$, the regularized-billiard Hamiltonian system \eqref{e:qham},
% \[
%\partial_z\begin{pmatrix}\u\\ \v\end{pmatrix}=\begin{pmatrix}\v\\ \nabla_\u B_\delta(\u)\end{pmatrix},
%\]
admits a locally unique orbit $\Psi_{h,\delta}=(\psi_{h,\delta},\partial_z\psi_{h,\delta})^{\rm T}$ homoclinic to zero for which
\beq
\lim_{\delta\rightarrow0+}\|\Psi_{h,\delta}(z)-\Psi_h(z)\|_{L^\infty}=0.
\eeq
\end{Proposition}

To prove the proposition, we recall  a Lemma, which states that restricted to inner points, the smooth regularized-billiard flow, denoted as ${\bf f}_\delta(z,\U)$, converges to the 
Birkhoff-billiard flow, denoted by ${\bf f}(z,\U)$, in the $C^\ell$-topology, as $\delta\rightarrow 0+$. We use the notation ${\bf f}_0(z,\U)={\bf f}(z,\U)$ where convenient.
\begin{Lemma}[{\bf \cite{turaev_1998}}]\label{l:bill}
For any inner point $\U_1$ and any $z>0$ so that $\U_2:={\bf f}(z,\U_1)$ is also an inner point,
there exists a neighborhood of $\U_{1\backslash 2}$, denoted as $\scn(\U_{1\backslash 2})$, and an interval containing $z$, denoted as $I$, so that, for sufficiently small $\delta>0$,
\beq
\|{\bf f}_\delta-{\bf f}\|_{C^\ell(I\times\scn(\U_1), \text{ } \scn(\U_2))}=\caO(\delta).
\eeq
\end{Lemma}

 We refer to \cite{turaev_1998} for the proof of Lemma \ref{l:bill} and give the proof of Proposition \ref{p:per} as follows.
\begin{Proof}
The idea is to show that the transversality stated in the Assumption \ref{a:trans} still holds for the regularized-billiard system \eqref{e:qham} for $\delta$ small. 
We only prove the case when $n=2$. The general proof for $n\in\Z^+$ is quite similar and thus omitted. 

The unperturbed ODE \eqref{e:ham} is Hamiltonian and thus admits $\mathscr{H}_0=\{\u\mid H(\u)=0\}$ as a three dimensional invariant manifold. For small $\xi>0$, the set
 \[\mathscr{S}(\xi):=\{\U\in\mathscr{H}_0\mid \u\in\sca, |\U-\Psi_{h,+}(z_2)|<\xi\},\]
where $\Psi_{h,+}(z_2):=\lim_{z\rightarrow z_2+}\Psi_h(z)$, provides a two dimensional cross section of $\Psi_h$ in $\mathscr{H}_0$.  However, 
the set $\mathscr{S}(\xi)$ has a drawback: As we perturb $\delta$ away from $0$, the intersections $\mathscr{S}(\xi)\cap \mathscr{M}^{s\backslash u}(0;\delta)$ can be empty.
Instead, we take a cross section close to $\mathscr{S}(\xi)$, which accommodates the persistence of the intersection. 
We choose $ z_0>0$ so small that $\{\psi_h(z)\mid z\in[-z_0, z_2+z_0]\backslash\{0,z_2\}\}\subset\scb^+$
and the set
\[\mathscr{S}_2(\xi):=z_0(\partial_z\psi_h(z_2+z_0),0)^{\rm T}+\mathscr{S}(\xi),\]
is a translation of $\mathscr{S}(\xi)$ centered at $\Psi_h(z_2+z_0)$,
admitting the same transversal property as $\mathscr{S}(\xi)$. As a result, we take $\xi$ so small that $\mathscr{S}_2(\xi)\cap\mathscr{M}^{s\backslash u}(0)$ is a $C^\ell$-smooth curve, denoted as ${\bf G}_{s\backslash u}(y)$, $y\in Y\subset \R$. From Assumption \ref{a:trans} it is true that $\Psi_{h}(z_2+z_0)$ is the only intersection point of two curves ${\bf G}_{s}$ and ${\bf G}_{u}$. Without loss of generality, we assume ${\bf G}_{s}(0)={\bf G}_{u}(0)=\Psi_{h}(z_2+z_0)$. Most importantly, the intersection is transversal, or equivalently, ${\bf G}_{s}^\prime(0)$ and ${\bf G}_{u}^\prime(0)$ are linearly independent. Since $\{\psi_h(z)\mid z>z_2\}\cap \sca=\emptyset$, for sufficiently small $\delta$, we can parameterize $\mathscr{M}^s(0;\delta)\cap\mathscr{S}_2(\xi)$ by ${\bf G}_{s,\delta}(t)$, $t\in T$. Moreover, we have
\beq\label{e:sest}
\|{\bf G}_{s,\delta}-{\bf G}_s\|_{(C^\ell(T))^4}=\caO(\delta).
\eeq

In order to describe $\mathscr{M}^u(0;\delta)\cap\mathscr{S}_2(\xi)$, we introduce a family of \textit{Poincar\'{e} maps} to control the flow near the collision curve. Close to the first collision point, we take a two dimensional cross section centered at $\Psi_h(-z_0)$,
\[\mathscr{S}_1(\tilde{\xi}):=\{\U\in\mathscr{H}_0\mid \big[\U-\Psi_h(-z_0)\big] \cdot (\partial_z\psi_h(-z_0),0)^{\rm T}=0, |\U-\Psi_h(-z_0)|<\tilde{\xi}\},\]
where $\tilde{\xi}$ is sufficiently small so that the $C^\ell$-smooth Poincar\'{e} map
\begin{equation}
\begin{matrix}
\Pi:&\mathscr{S}_1(\tilde{\xi}) &\longrightarrow & \mathscr{S}_2(\xi)\\
&\U &\longmapsto & {\bf f}(Z(\U), \U),
\end{matrix}
\end{equation}
is well defined with $\Pi(\mathscr{S}_1(\tilde{\xi}))\subset(\mathscr{S}_2(\xi))^o$. Here $Z(\U)>0$ is the $C^\ell$-smooth time of first arrival at $\mathscr{S}_2(\xi)$. Based on Lemma \ref{l:bill} and the definition of $\Pi$, there exists $\tilde{\delta}_0\in(0,\delta_0)$ so that, for any $\delta\in[0,\tilde{\delta}_0]$, the Poincar\'{e} map
\begin{equation}
\begin{matrix}
\Pi_\delta:&\mathscr{S}_1(\tilde{\xi}) &\longrightarrow & \mathscr{S}_2(\xi)\\
&\U &\longmapsto & {\bf f}_\delta(Z_\delta(\U),\U),
\end{matrix}
\end{equation}
is well-defined with $\Pi_\delta(\mathscr{S}_1(\tilde{\xi}))\subset(\mathscr{S}_2(\xi))^o$. We now define $\widetilde{\bf G}_u(y):=\Pi^{-1}({\bf G}_u(y))$ and thus there exists an interval $\widetilde{Y}\subseteq Y$ such that the $C^\ell$-smooth curve $\widetilde{\bf G}_u(y)$, $y\in\widetilde{Y}$, is the intersection $\mathcal{W}^u(0)\cap\mathscr{S}_1(\tilde{\xi})$. Similarly to the case of $\mathcal{W}^s(0;\delta)$, for sufficiently small $\delta$, the intersection $\mathcal{W}^u(0;\delta)\cap\mathscr{S}_1$ can be parameterized by $\widetilde{\bf G}_{u,\delta}(y)$, $y\in\widetilde{Y}$. Moreover, we have
\beq\label{e:alpha}
\|\widetilde{\bf G}_{u,\delta}-\widetilde{\bf G}_u\|_{(C^\ell(\widetilde{Y}))^4}=\caO(\delta).
\eeq
It is now straightforward to see that ${\bf G}_{u,\delta}(y):=\Pi_\delta(\widetilde{\bf G}_{u,\delta}(y))$ is the parameterization of $\mathcal{W}^u(0;\delta)$ in $\Pi_\delta(\mathscr{S}_1(\tilde{\xi}))\subset\mathscr{S}_2(\xi)$. Furthermore,
we have
\beq\label{e:inq}
\begin{aligned}
\|{\bf G}_{u,\delta}-{\bf G}_u\|_{(C^\ell(\widetilde{Y}))^4}&=\|\Pi_\delta(\widetilde{\bf G}_{u,\delta})-\Pi(\widetilde{\bf G}_u)\|_{(C^\ell(\widetilde{Y}))^4}\\
&\leq\|\Pi_\delta(\widetilde{\bf G}_{u,\delta})-\Pi(\widetilde{\bf G}_{u,\delta})\|_{(C^\ell(\widetilde{Y}))^4}+\|\Pi(\widetilde{\bf G}_{u,\delta})-\Pi(\widetilde{\bf G}_u)\|_{(C^\ell(\widetilde{Y}))^4}.
\end{aligned}
\eeq
From Lemma \ref{l:bill} it is known that the first term on the right side in \eqref{e:inq} is of order $\caO(\delta)$ while the estimate \eqref{e:alpha} indicates that 
the second term on the right side in \eqref{e:inq} is also of order $\caO(\delta)$. As a result, we have
\beq\label{e:uest}
\|{\bf G}_{u,\delta}-{\bf G}_u\|_{(C^\ell(\widetilde{Y}))^4}=\caO(\delta).
\eeq
Combining the estimates \eqref{e:sest} and \eqref{e:uest}, we conclude that, for sufficiently small $\delta>0$, the manifolds $\mathscr{S}_1(\xi)\cap\mathscr{M}^{s\backslash u}(0,\delta)$ intersect transversally at a point. This intersection point corresponds to a homoclinic orbit, denoted as $\Psi_{h,\delta}$, converging to $\Psi_h$ as $\delta$ goes to zero.
 \end{Proof}

%%%%%%%%%%%%%%%%%%%%%%%%%%%%%%%%%%%%
\subsection{Spectral analysis of the linearization: the origin and the collisions}
%%%%%%%%%%%%%%%%%%%%%%%%%%%%%%%%%%%%
In this sub-section, we analyze the eigenvalue problem associated to the linearization of the regularized Birkhoff-Billiard homoclinics established in Proposition\,\ref{p:per}. More specifically, for $\psi_{h,\delta}$ the homoclinic associated to the regularized Birkhoff-potential  $B_\delta$, then the corresponding linearization \eqref{e:L-delta} takes the form
\[
%\begin{matrix}
\mathcal{L}_\delta:=
%:&(H^2(\R))^2&\longrightarrow&(L^2(\R))^2\\&\w&\longmapsto&
\partial_z^2-\H_\delta(z),
%\end{matrix}
\]
where $\H_\delta(z):=\nabla_\u^2B_\delta(\psi_{h,\delta}(z))$ is the Hessian of $B_\delta(\u)$ with respect to $\u$ at $\psi_{h,\delta}(z)$.
The spectrum, $\sigma(\mathcal{L}_\delta)$, plays an important role in associated pearling stability of the corresponding quasi-bilayers $\u_q$ obtained by
dressing admissible interfaces with $\Psi_{h,\delta}(z)$.  The operator $\mathcal{L}_\delta$ is closed for $\delta>0$ but not well-defined for $\delta=0$. 
For small $\delta>0$, we localize the spectrum of the operator $\caL_\delta$, which consists of the essential spectrum $\sigma_{ess}(\caL_\delta)$ and the point spectrum $\sigma_{pt}(\caL_\delta)$.

%%%%%%%%%%%%%%%%%%%
%\paragraph{Essential spectrum and the origin}
%%%%%%%%%%%%%%%%%%%
By the Weyl essential spectrum theorem, see \cite{KapitulaPromislow_2013, kato, reedsimon78}, the essential spectrum of $\sigma_{ess}(\caL_\delta)$, is determined by the 
behavior of $B_\delta(\u)$ near the \textit{origin} $\u=0$. Since as $z\rightarrow\pm\infty$, the symmetric matrix $\H_\delta(z)$ limits to the constant matrix
$\nabla_\u^2B(0)=2{\bf I}_2$, independent of the choice of $\delta$, we introduce the \textit{asymptotic operator} associated with $\caL_\delta$,
\beq
\begin{matrix}
\mathcal{L}_\infty:&(H^2(\R))^2&\longrightarrow&(L^2(\R))^2\\
&\w&\longmapsto&\partial_z^2\w-2\w.
\end{matrix}
\eeq
It follows that $\caL_\delta-\caL_\infty$ is a compact operator and from the Weyl essential spectrum theorem the two operators share the same essential spectrum,
\beq\label{e:ess1}
\sigma_{ess}(\caL_\delta)=\sigma_{ess}(\caL_\infty)=(-\infty,-2].
\eeq

The point spectrum of a linear operator is in general more difficult to localize than its essential spectrum; see \cite{KapitulaPromislow_2013}. 
Instead of an exhaustive description of the point spectrum $\sigma_{\rm pt}(\caL_\delta)$, we characterize  the eigenvalues arising from \textit{collisions}, which
correspond to the largest eigenvalues in the limit $\delta\to0^+.$ 

We start by defining the collision times of perturbed homoclinics. For a perturbed homoclinic trajectory $\psi_{h,\delta}$ arising from the unperturbed 
$n$-collision homoclinic trajectory $\psi_h$ as in Proposition \ref{p:per}, the set of \textit{collision times} of the homoclinic $\psi_{h,\delta}$ are defined as 
\[
\{z\in\R\mid \partial_z\rho(\psi_{h,\delta}(z))=0,  \psi_{h,\delta}(z)\in\scn_\delta(\sca)\}.
\]
This set admits $n$ elements, denoted as $\{z_{i,\delta}\}_{i=1}^n$. For simplicity, we assume $z_{1,\delta}=0<\cdots<z_{n,\delta}$.
For a collision time $z_{i,\delta}$, the point ${\bf c}_{i,\delta}:=\psi_{h,\delta}(z_{i,\delta})$ is the \textit{collision point} and there is a \textit{collision interval} $  K_i(\delta)$ such that $z_{i,\delta}\in  K_i(\delta)$ and 
 \beq\label{e:cin}
 \bigcup_{i=1}^n  K_i(\delta)=\{z\in\R\mid \psi_{h,\delta}(z)\in\scn_\delta(\sca)\}.
 \eeq
By definition, we have $ \lim_{\delta\rightarrow0^+}z_{i,\delta}=z_i$ and $ \lim_{\delta\rightarrow0^+}{\bf c}_{i,\delta}={\bf c}_i$, where $\{z_i\}$ and $\{{\bf c}_i\}$  are collision times and points of the unperturbed homoclinic $\psi_h$.

\begin{figure}
  \begin{center}
    \includegraphics[width=0.88\textwidth]{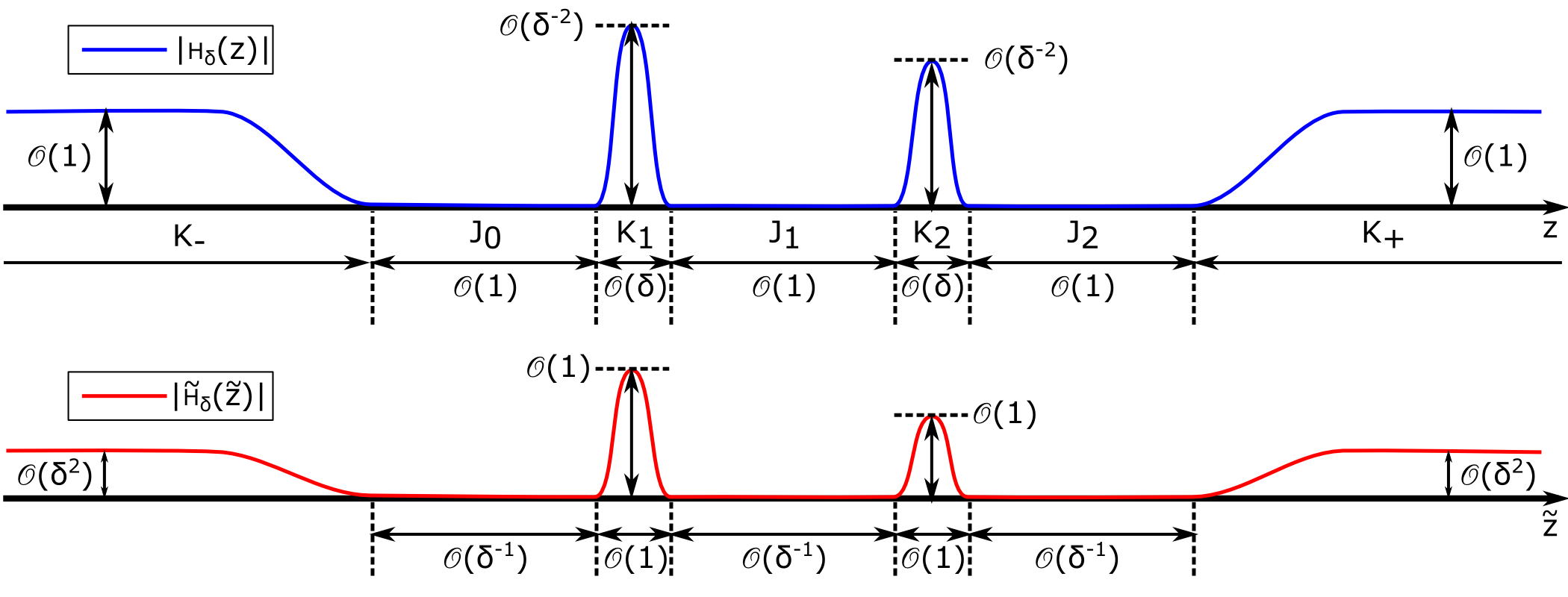}
  \end{center}
  \caption{The top and bottom panel are respectively cartoon sketches of the profiles $|\H_\delta |$ and $|\widetilde{\H}_\delta|$. }
\label{f:Hpro}
\end{figure}

To decompose the operator $\caL_\delta$ into a sum of collision operators and asymptotic operators, we chop the real line into $2n+3$ intervals; 
see the top panel of Figure \ref{f:Hpro} for an illustration of the case $n=2$. More specifically, we define
\begin{equation*}
\begin{cases}
Z_-(\delta):=\inf\{z\mid \psi_{h,\delta}(s)\in\scb^+, s\in[z, \min\{  K_1(\delta\}]\},\\
Z_+(\delta):=\sup\{z\mid \psi_{h,\delta}(s)\in\scb^+, s\in[\max\{  K_n(\delta)\}, z]\},\\
\end{cases}
\end{equation*}
and 
\[
K_-:=(-\infty, Z_-],\text{ }
J_0:=(Z_-, \min  K_1), \text{ }
J_i:=(\max  K_i, \min  K_{i+1}),\text{ }
J_n:=(\max  K_n, Z_+), \text{ }
K_+:=[Z_+,+\infty),
\]
where $i=1,\cdots, n-1$.
As a result, we have 
\beq\label{e:decom}
\R=\left(\bigcup_{i=1}^n  K_i(\delta)\right) \cup \left(\bigcup_{i=\pm}K_i(\delta)\right) \cup \left(\bigcup_{i=0}^{n+1}J_i(\delta)\right) .
\eeq
Moreover, the length of every $J_i$ is of order $\caO(1)$. As for the length of $K_i$, by Proposition \ref{p:per}, we conclude that, for sufficiently small $\delta>0$,
$\int_{z\in  K_{1\backslash 2}(\delta)}|\partial_z\psi_{h,\delta}(z)|\rmd z=\caO(\delta)$ and $\min_{z\in  K_{1\backslash 2}(\delta)}|\partial_z\psi_{h,\delta}(z)|=\caO(1)$,
yielding 
\beq
|  K_i|=\caO(\delta), i=1,\cdots, n.
\eeq

Let $\{\kappa_i(z)\}_{i=\pm, 1}^n$ be a partition of unity of the real line such that 
\beq
\begin{cases}
\sum_{i=\pm, 1}^n \kappa_i(z)=1, &z\in\R,\\
 \kappa_i(z)=1, &z\in K_i(\delta), i=1, \cdots, n, \pm.
\end{cases}
\eeq
Using such a partition of unity, we define 
\beq
\caK_{i,\delta}:=\kappa_i\caL_\delta, \quad  \caA_{\pm,\delta}:=\kappa_\pm\caL_\delta.
\eeq
and call $\caK_{i,\delta}$ the $i$-th \textit{collision operator}  and $\caA_{\pm,\delta}$ the \textit{positive$\backslash$negative asymptotic operator.}
 As a result, we have the decomposition
\beq
\caL_\delta=\caA_{-,\delta}+ \caA_{+,\delta}+\sum_{i=1}^n\caK_{i,\delta}.
\eeq

In the remainder of this section, we show that for small $\delta>0$, each collision operator admits a positive eigenvalue which persists in the whole linearized operator $\caL_\delta$. 
To start, we take a close look at $\H_\delta(z)$. A straightforward calculation shows that
\beq\label{e:H}
\H_\delta(z)=\begin{cases}
H_1(z;\delta)[\nabla_\u\rho(\psi_{h,\delta}(z))]^{\rm T}\nabla_\u\rho(\psi_{h,\delta}(z))+H_2(z;\delta)\nabla_\u^2\rho(\psi_{h,\delta}(z)), & z\in   \bigcup_{i=1}^n K_i(\delta),\\
\nabla_\u^2W(\psi_{h,\delta}(z)), & z\in   K_-(\delta)\cup  K_+(\delta),\\
0, &z\in\bigcup_{i=0}^n  J_i(\delta),\\
\end{cases}
\eeq
where 
\[
H_1(z;\delta):=(b_-+b_+)h_\delta^\prime(\rho(\psi_{h,\delta}(z)))=\caO(\delta^{-2}), \quad H_2(z;\delta):=(b_-+b_+)h_\delta(\rho(\psi_{h,\delta}(z)))=\caO(\delta^{-1}).
\]
Since $\H_\delta=\caO(\delta^{-2})$ on the $\caO(\delta)$-long intervals $K_i(\delta)$, there is no obvious well-defined limits for 
$\caL_\delta$ and $\caK_{i,\delta}$ as $\delta$ goes to zero. Applying a rescaling $z=\delta \tilde{z}$ and a $\delta^2$ scaling to 
$\caL_\delta$, we define a rescaled linearized operator
\beq\label{e:rescale}
\begin{matrix}
\widetilde{\mathcal{L}}_\delta:&(H^2(\R))^2&\longrightarrow&(L^2(\R))^2\\
&\w(\tilde{z})&\longmapsto&\partial_{\tilde{z}}^2\w(\tilde{z})-\widetilde{\H}_\delta(\tilde{z})\w(\tilde{z}),
\end{matrix}
\eeq
where $\widetilde{\H}_\delta(\tilde{z}):=\delta^2\H_\delta(\delta \tilde{z})$; see the bottom panel in Figure \ref{f:Hpro} for a scalar description. 
The sets $\sigma_{pt}(\widetilde{\caL}_\delta)$ and $\sigma_{pt}(\caL_\delta)$ are related by the map
\beq
\sigma_{pt}(\caL_\delta)=\delta^{-2}\sigma_{pt}(\widetilde{\caL}_\delta).
\eeq
More specifically, given $\caL_\delta \u(z)=\lambda \u(z)$, we have $\widetilde{\caL}_\delta \u(\delta \tilde{z})=\delta^2\lambda \u(\delta \tilde{z})$. Similarly, we define rescaled collision operators and rescaled asymptotic operators, respectively, 
\[
\widetilde{\caK}_{i,\delta}(\tilde{z}):=\kappa_i(\delta \tilde{z})\widetilde{\caL}_\delta, \quad \widetilde{\caA}_{\pm,\delta}(\tilde{z}):=\kappa_\pm(\delta \tilde{z})\widetilde{\caL}_\delta,
\]
for which the equality $\widetilde{\caL}_\delta=\sum_{i=1}^n\widetilde{\caK}_{i,\delta}+\sum_{j=\pm}\widetilde{\caA}_{j,\delta}$ holds. 
More importantly, every rescaled collision operator $\widetilde{\caK}_{i,\delta}$, up to the scaled translation $\tilde{z}\mapsto\tilde{z}-\delta^{-1}z_{i,\delta}$, admits a 
well-defined limit as $\delta$ goes to zero. In fact, we define the \textit{limit collision operators} as follows.
\begin{Definition}\label{d:limcol}
The $i$-th limit collision operator, denoted as $\caK_i$, is defined as 
\beq
\begin{matrix}
\caK_i: &(H^2(\R))^2 & \longrightarrow & (L^2(\R))^2 \\
& \w(\tilde{z}) & \longmapsto & \big(\partial_{\tilde{z}}^2-(b_++b_-)h^\prime(Z_i|\tilde{z}|)\nabla_\u\rho({\bf c}_i)[\nabla_\u\rho({\bf c}_i)]^{\rm T}\big)\w(\tilde{z}),
\end{matrix}
\eeq
where $Z_i:=\v_i\cdot\nabla_\u\rho({\bf c}_i)$,  $\v_i:=(\lim_{z\to z_i+}\partial_z\psi_h(z)-\lim_{z\to z_i-}\partial_z\psi_h(z))/2$, and $i=1,\cdots, n$.
\end{Definition}

We now have the following lemmas.
\begin{Lemma}\label{l:collision}
For small $\delta>0$ and any $i=1, \cdots, n$, we have 
\[
\lim_{\delta\to0+}\|\mathcal{T}_{i,\delta}\circ\tilde{\caK}_{i,\delta}-\caK_i(\tilde{z})\|_{B([H^2(\R)]^2,[L^2(\R)]^2)}=0,
\]
where $(\mathcal{T}_{i,\delta}\w)(\tilde{z}):=\w(\tilde{z}+\delta^{-1}z_{i,\delta})$ for any $\w\in(L^2(\R))^2$.
Moreover, the limit collision operator $\caK_i$ admits a unique positive eigenvalue, denoted $\nu_i$, that is,
\[
\sigma(\caK_i)\cap\{\nu\in\C\mid \Re\nu>0\}=\{\nu_i\}.
\]
Meanwhile, for small $\delta>0$, there exists $0<\nu_{i,\delta}\in\sigma_{pt}(\widetilde{\caK}_{i,\delta})$ such that $\nu_i=\lim_{\delta\to 0+}\nu_{i,\delta}$ and 
 $\nu_{i,\delta}$ is the unique positive eigenvalue of $\widetilde{\caK}_{i,\delta}$ with order $\caO(1)$.
\end{Lemma}

\begin{Proof}
We prove the result for the scaled first-collision operator and omit the rest. The truncated first-collision Hessian $\widetilde{\H}_{1,\delta}(\tilde{z}):=\kappa_{1,\delta}(\delta\tilde{z})\widetilde{\H}_\delta(\tilde{z})$ admits the expression
\beq\label{e:H1}
\widetilde{\H}_{1,\delta}(\tilde{z})=
\begin{cases}
\delta^2\H_\delta(\delta\tilde{z}),  &\tilde{z}\in \delta^{-1}\caK_{1,\delta},\\
0, &\tilde{z}\notin \delta^{-1}\caK_{1,\delta},\\
\end{cases}
\eeq
Applying the rescalings $\rho=\delta\tilde{\rho}$ and $ \u= {\bf c}_{1,\delta}+\delta\tilde{\u}$ to the rescaled level-set function $\tilde{\rho}(\tilde{\u}):=\delta^{-1}\rho({\bf c}_{1,\delta}+\delta\tilde{\u})$,
we have that, for $\widetilde{\psi}_{h,\delta}(\tilde{z}):=\delta^{-1}(\psi_{h,\delta}(\delta\tilde{z})-{\bf c}_{1,\delta})$ and $\tilde{z}\in \delta^{-1}\caK_{1,\delta}$,
\begin{equation}
\label{e:limit}
\begin{cases}
\nabla_\u\rho(\psi_{h,\delta}(z))=\nabla_{\tilde{\u}}\tilde{\rho}(\widetilde{\psi}_{h,\delta}(\tilde{z}))=\nabla_\u\rho({\bf c}_1)+\caO(\delta),\\
\tilde{\rho}(\widetilde{\psi}_{h,\delta}(\tilde{z}))=Z_1|\tilde{z}|+\caO(\delta),\\
\delta^2h_\delta^\prime(\rho(\psi_{h,\delta}(z)))=h^\prime(\tilde{\rho}(\widetilde{\psi}_{h,\delta}(\tilde{z})))=h^\prime(Z_1|\tilde{z}|)+\caO(\delta). \\
\end{cases}
\end{equation}
 Inserting the expansion \eqref{e:limit} and the expression \eqref{e:H} into \eqref{e:H1}, we obtain
\[
\lim_{\delta\to 0+}\widetilde{\H}_{1,\delta}(\tilde{z})=(b_++b_-)h^\prime(Z_1|\tilde{z}|)\nabla_\u\rho({\bf c}_1)[\nabla_\u\rho({\bf c}_1)]^{\rm T},
\]
implying that $\lim_{\delta\to0+}\caK_{1,\delta}=\caK_1$.

We are left to show that $\caK_1$ admits a positive eigenvalue. Defining 
\[
{\bf n}_1:=\frac{\nabla_\u\rho({\bf c}_1)}{\|\nabla_\u\rho({\bf c}_1)\|}, \quad {\bf t}_1:=\frac{\lim_{z\to0+}\partial_z\psi_h(z)+\lim_{z\to0-}\partial_z\psi_h(z)}{\|\lim_{r\rightarrow0+}\partial_z\psi_h(z)+\lim_{r\rightarrow0-}\partial_z\psi_h(z)\|}, \quad \mathcal{N}_1:=\begin{pmatrix}{\bf n}_1&{\bf t}_1\end{pmatrix},
\]
we apply the change of coordinates $\w=\mathcal{N}_1\widetilde{\w}$ to the operator $\widetilde{\caK}_{1,\delta}$ to find
\[
\cad_{1,\delta}:=\mathcal{N}_1^{-1}\circ\widetilde{\caK}_{1,\delta}\circ \mathcal{N}_1=
\begin{pmatrix}
\partial_{\tilde{z}}^2-(b_++b_-)h^\prime(Z_1\tilde{z}) & 0\\
0 & \partial_{\tilde{z}}^2
\end{pmatrix}.
\]
We note that 
\[
\big( \partial_{\tilde{z}}^2-(b_++b_-)h^\prime(Z_1\tilde{z})  \big) \widetilde{w}_0(\tilde{z})=0,
\]
where $\widetilde{w}_0(\tilde{z}):=\lim_{\delta\to0+}  {\bf n}_1\cdot \partial_z\psi_{h,\delta}(\delta \tilde{z})$ is an increasing function with $\widetilde{w}_0(0)=0$. By Sturm-Louville theory, the operator $\partial_{\tilde{z}}^2-(b_++b_-)h^\prime(Z_1\tilde{z}) $, and thus, $\caK_1$, admits one and only one positive eigenvalue $\nu_1$. The continuation of $\nu_1$ in $\widetilde{\caK}_{1,\delta}$ is a simple classical perturbation argument \cite{kato}, which concludes the proof.
\end{Proof}

\begin{Lemma}\label{l:asymp}
For small $\delta>0$, the asymptotic operators $\widetilde{\caA}_{\pm,\delta}$ have no $\caO(1)$ positive eigenvalues,
\[
\sigma_{pt}(\widetilde{\caA}_{\pm,\delta})\cap \{\nu\in\R^+\mid \nu=\caO(1)\}=\emptyset.
\]
\end{Lemma}

\begin{Proof}
We argue by contradiction to show that
\[
\sigma_{pt}(\caA_{-,\delta})\cap \{\nu\in\R^+\mid \nu=\caO(\delta^{-2})\}=\emptyset.
\]
 The proof for the operator $\caA_{+,\delta}$ is essentially the same and thus omitted.
 Due to the radial symmetry of $B(\u)$ near the origin,  there exist $\alpha(\delta)\in[0,\pi/2]$ so that the closure of the trajectory $\{\psi_{h,\delta}(z)\mid z\in K_-\}$ is a line segment through the origin whose angle with respect to $u_1$-axis in the $\u$-plane is $\alpha$. Consequently, denoting ${\bf t}_-:=(\cos(\alpha),\sin(\alpha))^{\rm T}$,
 $r(z;\delta):=|\psi_{h,\delta}(z)|$ and $Z_m(\delta):=\{z\in K_-\mid r(z)=R_0/2\}$, we have
 \begin{equation} 
 \H_{-,\delta}(z)=
 \begin{cases}
 2{\bf I}_2, & z\in(-\infty, Z_m),\\
 (b^{\prime\prime}(r(z))-b^{\prime}(r(z))/r(z)){\bf t}_-{\bf t}_-^{\rm T}+b^{\prime}(r(z))/r(z){\bf I}_2, & z\in[Z_m, Z_-],\\
0, &z\in[Z_-,\infty).
 \end{cases}
 \end{equation}
We denote the $2\times 2$ rotation matrix with angle $\alpha$ as $\mathcal{N}_-$ and apply the change of coordinates $\w=\mathcal{N}_-\widetilde{\w}$ to the operator $\caA_{-,\delta}$, yielding
\[
\cad_{-,\delta}:=\mathcal{N}_-^{-1}\circ\caA_{-,\delta}\circ \mathcal{N}_-=
\begin{pmatrix}
\partial_z^2-b_1(z) & 0\\
0 & \partial_z^2-b_2(z)
\end{pmatrix},
\]
where 
  \begin{equation*}
  b_1(z)=
 \begin{cases}
 2, & z\in(-\infty, Z_m),\\
 b^{\prime\prime}(r(z)), & z\in[Z_m, Z_-],\\
0, &z\in[Z_-,\infty),
 \end{cases}
 \quad  
 b_2(z)=
 \begin{cases}
2, & z\in(-\infty, Z_m),\\
b^{\prime}(r(z))/r(z), & z\in[Z_m, Z_-],\\
0, &z\in[Z_-,\infty).
 \end{cases}
 \end{equation*}
Assume that there exists a positive constant $\nu=\caO(1)$ so that $\delta^{-2}\nu\in \sigma_{pt}(\partial_z^2-b_1(z))$ with an eigenfunction $u(z)\in L^2(\R)$, that is,
$\partial_z^2u(z)-b_1(z)u(z)=\delta^{-2}\nu u(z)$. It is then straightforward to see that the \textit{angle function},
$\beta(z):=\frac{u(z)}{\partial_zu(z)}$, solves the following second order ODE problem
\begin{subnumcases}{\label{e:angleode} }
 \partial_z\beta=(\delta^{-2}\nu+b_1(z))-\beta^2, &\label{e:ode}\\
 \beta(Z_m)=\sqrt{\delta^{-2}\nu+2}, \beta(Z_-)=-\delta^{-1}\sqrt{\nu}, &\label{e:intl}
\end{subnumcases}
which is impossible for sufficiently small $\delta$. The reason is quite simple: since $b_1=\caO(1)$ for all $z\in\R$, we have $\partial_z\beta\gg 1$ when $\beta=0$, 
which implies that any forward solution of the ODE \eqref{e:ode} with positive initial condition will stay positive for all forward time. 
A similar argument also applies to the operator $\partial_z^2-b_2(z)$. Therefore, the operator $\cad_{-,\delta}$, thus $\caA_{-,\delta}$, does not 
admit any positive eigenvalue of order $\caO(\delta^{-2})$, which concludes the proof.
\end{Proof}

Having established Lemmas \ref{l:collision} and \ref{l:asymp}, we are ready to prove the main theorem of this section.
\begin{Proof}[Proof of Theorem \ref{t:eigen}] 
It remains to prove the statements about the collision eigenvalues. 
%We first claim that the asymptotic operators $\caA_{\pm,\delta}$ has no impact on collision eigenvalues. In other words,  $\lambda$ is a collision eigenvalue of the operator $\caL_\delta$ if and only if $\lambda$ is a collision eigenvalue of the \textit{total collision} operator $\caK_\delta:=\sum_{j=1}^n\caK_{j,\delta}$. The ``only if" direction can be readily shown by an argument similar to the one about the angle function $\beta$ in the proof of Lemma \ref{l:asymp}. 
%--the rescaled \textit{total collision} operator $\widetilde{\caK}_\delta:=\sum_{j=1}^n\widetilde{\caK}_{j,\delta}$. Without loss of generality, we only prove the case $n=2$. The case $n>2$ can be readily checked via a straightforward inductive argument.
%
Without loss of generality, we consider only the case $n=2$. Recall that the rescaled linearized operator $\widetilde{\caL}_\delta$ admits the decomposition
\[
\widetilde{\caL}_\delta=\widetilde{\caA}_{-,\delta}+\widetilde{\caK}_{1,\delta}+\widetilde{\caK}_{2,\delta}+\widetilde{\caA}_{+,\delta},
\]
where the rescaled collision operator $\widetilde{\caK}_{i,\delta}$ admits a positive eigenvalue $\nu_{i,\delta}=\caO(1)$, as shown in Lemma \ref{l:collision}. We restrict ourselves to the continuation of $\nu_{1,\delta}$ in $\widetilde{\caL}_\delta$. The argument is inductive: we first show  the continuation of $\nu_{1,\delta}$ in  $\widetilde{\caA}_{-,\delta}+\widetilde{\caK}_{1,\delta}$, then similar arguments can apply to $\widetilde{\caA}_{-,\delta}+\widetilde{\caK}_{1,\delta}+\widetilde{\caK}_{2,\delta}$ and eventually $\widetilde{\caL}_\delta$. Indeed, 
the proof boils down to showing the existence of an eigenvalue $\widetilde{\lambda}_{1,\delta}\in \sigma_{pt}(\widetilde{\caK}_{1,\delta}+\widetilde{\caA}_{-,\delta})$ close to $\nu_{1,\delta}$,

To illustrate the ideas without getting involved in technicalities, we simplify the argument by assuming that $\widetilde{\caA}_{-,\delta}$ and $\widetilde{\caK}_{1,\delta}$ are 
scalar operators.  Defining $\kappa_1:=\min K_{1,\delta}$, for $\tilde{\nu}$ small, we aim to locate initial conditions $(d_1(\tilde{\nu}),d_2(\tilde{\nu}))$ such that 
the solution $(\w(\tilde{z},d_1,d_2), \partial_{\tilde{z}}\w(\tilde{z}, d_1,d_2))$ of the initial value problem
\begin{equation*}
\begin{cases}
\widetilde{\caK}_{1,\delta}\w=(\nu_{1,\delta}+\tilde{\nu})\w,\\
\w(\delta^{-1}\kappa_1)=d_1,\\
\partial_{\tilde{z}}\w(\delta^{-1}\kappa_1)=d_2,\\
\end{cases}
\end{equation*}
satisfies
\[
\lim_{\tilde{z}\to +\infty}(\w(\tilde{z},d_1,d_2), \partial_{\tilde{z}}\w(\tilde{z}, d_1,d_2))=0.
\]
Since $\nu_{1,\delta}\in\sigma_{pt}(\widetilde{\caK}_{1,\delta})$, we define 
\[
(d_1(0),d_2(0))=(1,\sqrt{\nu_{1,\delta}}),
\]
For $\tilde{\nu}$ small, there exists a scalar function $f(\tilde{\nu})=f_0\tilde{\nu}+\caO(\tilde{\nu}^2)$, $f_0\neq0$, such that 
\[
(d_1(\tilde{\nu}),d_2(\tilde{\nu}))=(1,\sqrt{\nu_{1,\delta}+\tilde{\nu}})+f(\tilde{\nu})(1,-\sqrt{\nu_{1,\delta}+\tilde{\nu}}),
\]
The fact that $f_0\neq0$ is a straightforward argument based on Pr\"{u}fer's substitution \cite{BirkhoffRota_89} and the comparison theorem.
Therefore, we have 
\begin{equation}\label{e:jump}
\begin{aligned}
&(\w(\delta^{-1}R_-,d_1(\tilde{\nu}),d_2(\tilde{\nu})), \partial_{\tilde{z}}\w(\delta^{-1}R_-, d_1(\tilde{\nu}),d_2(\tilde{\nu})))\\
=&\rme^{\delta^{-1}\sqrt{\nu_{1,\delta}+\tilde{\nu}}(R_--\kappa_1)}(1,\sqrt{\nu_{1,\delta}+\tilde{\nu}})+\rme^{\delta^{-1}\sqrt{\nu_{1,\delta}+\tilde{\nu}}(\kappa_1-R_-)}f(\tilde{\nu})(1,-\sqrt{\nu_{1,\delta}+\tilde{\nu}}).
\end{aligned}
\end{equation}
Noting that $\nu_{1,\delta}=\caO(1)$, we have
\[
\sqrt{\nu_{1,\delta}}/\delta\gg 1.
\]
Therefore, for small $\delta>0$, there always exist $\tilde{\nu}_1$ such that the solution to the initial value problem
\begin{equation*}
\begin{cases}
(\widetilde{\caK}_{1,\delta}+\caA_{-,\delta})\v=(\nu_{1,\delta}+\tilde{\nu})\v,\\
(\v,\partial_{\tilde{z}}\v)(\delta^{-1}R_-,d_1(\tilde{\nu}),d_2(\tilde{\nu}))=(\w(\delta^{-1}R_-,d_1(\tilde{\nu}),d_2(\tilde{\nu})), \partial_{\tilde{z}}\w(\delta^{-1}R_-, d_1(\tilde{\nu}),d_2(\tilde{\nu}))),
\end{cases}
\end{equation*}
satisfying
\[
\lim_{\tilde{z}\to -\infty}(\v(\tilde{z},d_1(\tilde{\nu}_1),d_2(\tilde{\nu}_1)), \partial_{\tilde{z}}\v(\tilde{z}, d_1(\tilde{\nu}_1),d_2(\tilde{\nu}_1)))=0.
\]
This is precisely what is required for $\tilde{\lambda}_{1,\delta}=\nu_{1,\delta}+\tilde{nu}_1$ to be an eigenvalue of $\widetilde{\caK}_{1,\delta}+\caA_{-,\delta}$, and thus concludes the proof.
\end{Proof}
%Given 
%\[
%\widetilde{  K}_1:=\{s\in[s_{1,-},s_{1,+}]\mid s_{1,-}=\lim_{\delta\rightarrow 0+}\delta^{-1}\min\{r\in  K_1\}, s_{1,+}=\lim_{\delta\rightarrow 0+}\delta^{-1}\max\{r\in  K_1\}\},
%\]
%we define
%\beq
%\widetilde{\H}(s):=\lim_{\delta\rightarrow0+}\widetilde{\H}_\delta(s)=
%\begin{cases}
%\widetilde{\Psi}_1(s), & s\in \widetilde{  K}_1,\\
%0, & s\notin \widetilde{  K}_1,
%\end{cases}
%\eeq
%where 
%\[
%\widetilde{\Psi}_1(s):=\lim_{\delta\rightarrow0+}\delta^2\Psi_1(\delta s; \delta).
%\]
%\begin{Lemma}
%for
%\end{Lemma}

%%%%%%%%%%%%%%%%%%%%%%%%%%%%%%%%%%%%%%%%%%%%%%%%%%%%%%
\section{Geometric evolution of multi-component bilayers}\label{s:4}
\begin{figure}[!ht]
  \begin{center}
    \includegraphics[width=0.8\textwidth]{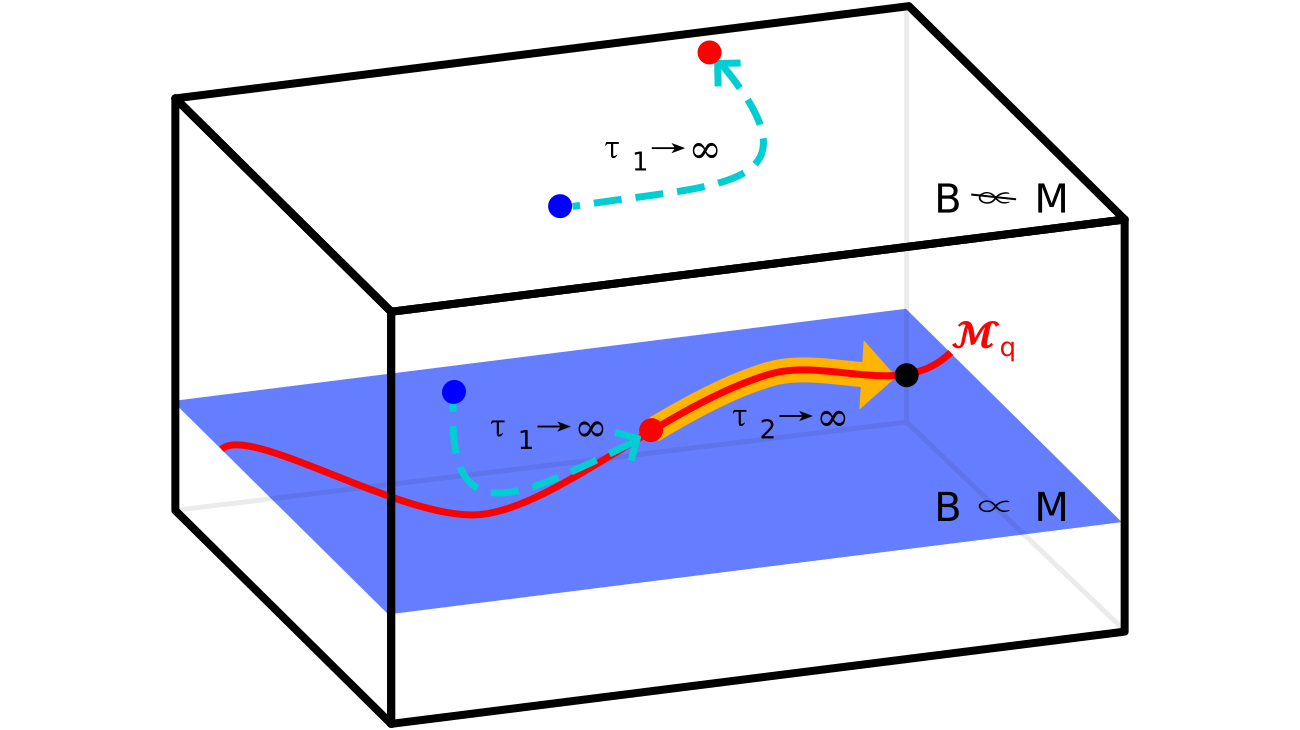}
  \end{center}
  \caption{   A reduced dimension ``cartoon'' depiction of the multi-time-scale geometric flow that drives pseudo-bilayers onto the affine subspace of quasi-bilayers.  
  The cuboid represents the set of pseudo-bilayers, the blue surface represents the set of pseudo-bilayers with leading order background state $\veps {\bf B}$ proportional to the quasi-bilayer mass vector ${\bf M}$ and the red curve represents the set of quasi-bilayers, $\mathscr{M}_q$. On the fast time scale, $\tau=\veps t$,
 the light-blue dashed curve depicts the evolution of a typical pseudo-bilayer initial data with $ {\bf B}\propto {\bf M}$  (deep blue dot on the blue surface), which converges to its fast time equilibrium, (red dot); pseudo-bilayers with $\veps {\bf B}$ not proportional to ${\bf M}$ cannot converge to slow time equilibrium not on $\mathscr{M}_q$. While $\mathscr{M}_q$
 is composed of fast-time equilibria,  large transients occur within $\mathscr{M}_q$ on the slow time scale $\tau_2=\veps^2 t$, governed by the mass preserving Willmore flow,
 (\ref{e:nv2}),  depicted by thick orange curve connecting the fast time equilibrium (red dot) to its slow time equilibrium (black dot).  
}
\label{f:geomnfd}
\end{figure}

%%%%%%%%%%%%%%%%%%%%%%%%%%%%%%%%%%%%%%%%%%%%%%%%%%%%%%
In this section we analyze the geometric evolution of co-dimensional one quasi-bilayers, deriving their curvature-driven motion under the
$H^{-1}$ gradient flow of the weak mFCH. The results are formal in the sense that we assume the stability of the underlying quasi-bilayers, in particular their
stability to the pearling bifurcations described in Section\,\ref{s:3}.  Moreover, as is typical in a multi-scale analysis, we assume that the evolution
occurs at distinct time scales, which scale according to,
\[
\tau_\alpha=\veps^\alpha t, \quad \alpha=-2,-1,0,1,2.
\] 
More specifically, given a fixed time scale $\tau_\alpha=\veps^\alpha t$,
we adopt the asymptotic analysis techniques from \cite{DaiPromislow_2013, Pego_1989} to analyze the initial value problem of the rescaled weak mFCH equation,
\beq\label{e:intlh-1}
%\begin{cases}
\veps^\alpha\u_{\tau_\alpha}(x,\tau_\alpha)=\Delta \mu(x,\tau_\alpha),\\
%\u(x,0)=\u_q(x,\veps; {\bf m}, \Gamma),
%\end{cases}
\eeq
subject to initial data that corresponds to the equilibrium of the preceding time-scale. Here the chemical potential $\mu=\frac{\delta \mathcal{F}_M}{\delta \u}$ and the initial profile is a quasi-bilayer $\u_q(x,\veps; {\bf m}, \Gamma)$ with background state ${\bf m}$ and interface $\Gamma$, as defined in \eqref{e:quasib}. Since quasi-bilayers
are leading order equilibria of (\ref{e:intlh-1}) for $\alpha\leq 1$, we first consider the time scale $\tau_1=\veps t$, under a slightly broader 
class of initial data, called pseudo-bilayers and denoted as $\u_p(x,\veps; {\bf B}, \Gamma)$,
taking the form
\beq\label{e:pseudob}
\u_p(x, \veps; {\bf B},\Gamma)=\begin{cases}
\phi_h(z(x))+\veps(\caL_{0,\perp})^{-2}({\bf B}-\caL_0{\bf V}(\phi_h))(z(x))+\caO(\veps^2), & x\in\Gamma_{l_0},\\
\veps {\bf B}+\caO(\veps^2), & x\in\Omega\backslash \Gamma_{3l_0},
\end{cases}
\eeq
with a smooth transition in  $\Gamma_{3l_0}\backslash \Gamma_{l_0}$ as in \eqref{e:quasib}. Here $\phi_h$ is an orbit of \eqref{e:leadode} which is homoclinic to the origin,
and the pseudo-bilayer reduces to a quasi-bilayer when ${\bf B}=0$ to leading order.  We show that any pseudo-bilayer converges to a quasi-bilayer at leading order 
as $\tau_1$ goes to infinity, if and only if the leading order background state ${\bf B}$ is proportional to the quasi-bilayer mass vector ${\bf M}=\int_\R \phi_h(z)\rmd z$, 
see Figure \ref{f:geomnfd}. This condition represents a mass constraint -- the qausi-bilayer can be the dominant repository of mass only if the initial data has \textcolor{blue}{a} compatible mass constitution. On the longer time-scale $\tau_2=\veps^2 t$, we show that any quasi-bilayer profile initial data evolves as a quasi-bilayer at leading order 
subject to a mass preserving Willmore flow that incorporates the intrinsic curvature of the associated quasi-bilayer $\Phi_h$, defined in (\ref{e:aexp}).

In the outer region $\Omega\backslash \Gamma_{l_0}$  we have the expansions,
\begin{equation*}
\begin{cases}
\u(x,\tau_\alpha;\veps)&=\u_0+\veps\u_1+\veps^2\u_2+\veps^3\u_3+\caO(\veps^4),\\
\mu(x,\tau_\alpha;\veps)&=\mu_0+\veps\mu_1+\veps^2\mu_2+\veps^3\mu_3+\caO(\veps^4).
\end{cases}
\end{equation*}
while in the inner region $\Gamma_{l_0}$ the expansions take the form
\begin{equation*}
\begin{cases}
\u(x,\tau_\alpha;\veps)&=\tilde{\u}(z,s,\tau_\alpha;\veps)=\tilde{\u}_0+\veps\tilde{\u}_1+\veps^2\tilde{\u}_2+\veps^3\tilde{\u}_3+\caO(\veps^4),\\
\mu(x,\tau_\alpha;\veps)&=\tilde{\mu}(z,s,\tau_\alpha;\veps)=\tilde{\mu}_0+\veps\tilde{\mu}_1+\veps^2\tilde{\mu}_2+\veps^3\tilde{\mu}_3+\caO(\veps^4).
\end{cases}
\end{equation*}
The time derivative involving the normal velocity of the interface admits the expansion,
\[
\veps^\alpha\u_{\tau_\alpha}(x,\tau_\alpha;\veps)=\veps^q\left[\partial_{\tau_\alpha} \tilde{\u}+\partial_s\tilde{\u}\partial_{\tau_\alpha} s+\partial_z\tilde{\u}\partial_{\tau_\alpha} r\right]=-\veps^{q-1} V_{n,0} \partial_z\tilde{\u}+\caO(\veps^q),
\]
where $\mathbf{n}\cdot\partial_{\tau_\alpha}\Xi:=V_n(s,\tau_\alpha;\veps)=V_{n,0}+\veps V_{n,1}+\caO(\veps^2)$.
More specifically, we have
\begin{equation*}
\begin{cases}
\mu_0=&\nabla_\u^2W(\u_0)\nabla_\u W(\u_0),\\
\mu_1=&-\nabla_\u^2W(\u_0)\left( -\nabla_\u^2W(\u_0)\u_1+{\bf V}(\u_0)\right)-\left[-\nabla_\u^3W(\u_0)\u_1+(\nabla_u{\bf V}(\u_0))^{\rm T}\right]\nabla_\u W(\u_0),\\
\mu_2=&-\nabla_\u^2W(\u_0)\left[ \Delta\u_0-\nabla_\u^2W(\u_0)\u_2-\frac{1}{2}\nabla_\u^3W(\u_0)(\u_1,\u_1)+\nabla_\u{\bf V}(\u_0)\u_1\right]\\
&+\left[-\nabla_\u^3W(\u_0)\u_1+(\nabla_u{\bf V}(\u_0))^{\rm T}\right]\left( -\nabla_\u^2W(\u_0)\u_1+{\bf V}(\u_0)\right)\\
&-\left[\Delta-\nabla_\u^3W(\u_0)\u_2-\frac{1}{2}\nabla_\u^4W(\u_0)(\u_1,\u_1)+(\nabla_\u^2{\bf V}(\u_0)\u_1)^{\rm T}\right]\nabla_\u W(\u_0)-\eta_2\nabla_\u W(\u_0),
\end{cases}
\end{equation*}
and 
\begin{equation*}
\begin{cases}
\tilde{\mu}_0=&\Rmnum{1}_0\Rmnum{2}_0,\\
\tilde{\mu}_1=&\Rmnum{1}_0\Rmnum{2}_1+\Rmnum{1}_1\Rmnum{2}_0,\\
\tilde{\mu}_2=&\Rmnum{1}_0\Rmnum{2}_2+\Rmnum{1}_1\Rmnum{2}_1+\Rmnum{1}_2\Rmnum{2}_0+\Rmnum{3}_2,\\
\tilde{\mu}_3=&\Rmnum{1}_0\Rmnum{2}_3+\Rmnum{1}_1\Rmnum{2}_2+\Rmnum{1}_2\Rmnum{2}_1+\Rmnum{1}_3\Rmnum{2}_0+\Rmnum{3}_3,\\
\end{cases}
\end{equation*}
where 
\begin{equation*}
\begin{cases}
\partial_z^2+\veps H(z,s)\partial_z+\veps^2\Delta_G-\nabla_\u^2W(\tilde{\u})+\veps(\nabla_\u {\bf V}(\tilde{\u}))^{\rm T}:=\Rmnum{1}=\Rmnum{1}_0+\veps\Rmnum{1}_1+\cdots,\\
\left(\partial_z^2+\veps H(z,s)\partial_z+\veps^2\Delta_G\right) \tilde{\u}-\nabla_\u W(\tilde{\u})+\veps {\bf V}(\tilde{\u}):=\Rmnum{2}=\Rmnum{2}_0+\veps\Rmnum{2}_1+\cdots,\\
-\veps^2\left[ -\eta_1\left(\partial_z^2+\veps H(z,s)\partial_z+\veps^2\Delta_G \right)\tilde{\u} +\eta_2\nabla_\u W(\tilde{\u}) \right]:=\Rmnum{3}=\Rmnum{3}_0+\veps\Rmnum{3}_1+\cdots,\\
\end{cases}
\end{equation*}
We follow the usual matching procedure, first employed for the Cahn-Hilliard equation in \cite{Pego_1989}, for the inner and outer expansions; that is, 
for any given $x\in\Gamma$, we require that
\[
(\u_0+\veps\u_1+\veps^2\u_2+\veps^3\u_3+\cdots)(x+\veps z\mathbf{n},\tau_\alpha)\approx(\tilde{\u}_0+\veps\tilde{\u}_1+\veps^2\tilde{\u}_2+\veps^3\tilde{\u}_3+\cdots)(z,s,\tau_\alpha).
\]
For the cases $\alpha=-2,-1,0, 1$, the calculation closely follows the scalar case presented in \cite{DaiPromislow_2013}, and shows that
any quasi-bilayer profile is an equilibrium of the corresponding leading-order evolution system. We provide details for
the slow time scale $\tau_1=\veps t$ with initial pseudo-bilayer profile and the slow time scale $\tau_2=\veps^2t$ with initial data a quasi-bilayer profile
obtained as equilibrium of the $\tau_1$ evolution.

%%%%%%%%%%%%%%%%%%%%%%
\paragraph{Time scale $\tau_1=\veps t$: a quenched mean curvature flow.}
%%%%%%%%%%%%%%%%%%%%%%
We study the initial-value problem, 
\beq\label{e:intlh-1p}
\begin{cases}
\veps\u_{\tau_1}(x,\tau_1)=\Delta \mu(x,\tau_1),\\
\u(x,0)=\u_p(x,\veps; {\bf B}_0, \Gamma_0).
\end{cases}
\eeq
 At leading order the outer expansion takes the form,
\begin{equation}
\begin{cases}
& \mu_0=0, \\
& \mu_1=\left(\nabla_\u^2W(0)\right)^2\u_1, \quad \Delta \mu_1=0.
\end{cases}
\end{equation}
while the leading order terms in the inner expansion  take the form,
\begin{equation}
\begin{cases}
& \tilde{\mu}_0=0, \\
&\partial_z^2 \tilde{\mu}_1=0,\\
&-V_{n,0}\partial_z\phi_h=H_0(s)\partial_z\tilde{\u}_1+\partial_z^2\tilde{\u}_2.
\end{cases}
\end{equation}
A analysis similar to that in \cite{DaiPromislow_2013} yields the following leading-order evolution system, for $\tau_1\geq0$,
\begin{subnumcases}{\label{e:q1} }
\tilde{\u}_0=\phi_h(z), \quad \u_0=0, &\label{e:q10}\\
 \mu_1=\tilde{\mu}_1=\mathbf{B}_1(\tau_1), &\label{e:q11}\\
 \tilde{\u}_1=(\caL_{0,\perp})^{-2}(\mathbf{B}_1(\tau_1)-\caL_0{\bf V}(\phi_h)),&\label{e:q12}\\
 \u_1=\left( \nabla_\u^2 W(0)\right)^{-2} \mathbf{B}_1(\tau_1),&\label{e:q13}\\
V_{n,0}=\frac{\mathbf{B}_1(\tau_1)\cdot \mathbf{M}}{M_2}H_0,&\label{e:q14}
\end{subnumcases}
where ${\bf B}_1(0)=(\nabla_\u W(0))^2{\bf B}_0$ and we recall that $\mathbf{M}=\int_\R \phi_h(z)\rmd z$ and $M_2=\int_\R |\phi_h(z)|^2\rmd z$.
From (\ref{e:q11}-\ref{e:q13}), we have
\[
\u(x,\tau_1;\veps)=\begin{cases}
\phi_h(z(x))+\veps\tilde{\u}_1(x,\tau_1)+\caO(\veps^2), &x\in \Gamma_{l_0},\\
\veps (\nabla_\u^2 W(0))^{-2}{\bf B}_1(\tau_1)+\caO(\veps^2), &x\in \Omega\backslash \Gamma_{l_0},
\end{cases}
\]
whose $\veps$-scaled mass vector $\mathcal{M}$, defined in (\ref{e:calM-def}), admits the expansion,
\[
\mathcal{M}(\tau_1)=|\Omega| (\nabla_\u^2 W(0))^{-2}{\bf B}_1(\tau_1)+|\Gamma|(\tau_1){\bf M}+\caO(\veps).
\]
Expanding the mass vector $\mathcal{M}=\mathcal{M}_1+\caO(\veps)$ and the interface area $|\Gamma|(\tau_1)=\gamma_0(\tau_1)+\caO(\veps)$, we have
\beq\label{e:M1}
\mathcal{M}_1=|\Omega| (\nabla_\u^2 W(0))^{-2}{\bf B}_1(\tau_1)+\gamma_0(\tau_1){\bf M}.
\eeq
 Under a prescribed normal velocity, the time derivative of the surface area equals the integral of the product of mean curvature and normal velocity over the surface, that is,
 \beq\label{e:surfevo}
 \partial_{\tau_1} |\Gamma|=\int_SH_0V_n\rmd s.
\eeq
At leading order, \eqref{e:surfevo} takes the form $\partial_{\tau_1} \gamma_0=M_2^{-1}({\bf B}_1(\tau_1)\cdot {\bf M})\int_SH_0^2\rmd s$. Defining the quantity
\[E(\tau_1):={\bf B}_1(\tau_1)\cdot {\bf M},\]
then (\ref{e:surfevo}) and \eqref{e:M1} prescribe the evolution of $E$, that is,
\beq\label{e:extra}
\partial_{\tau_1} E=-\frac{{\bf M}^{\rm T}\left[\nabla_\u^2 W(0)\right]^2{\bf M}}{|\Omega|M_2} E\int_S H_0^2\rmd s.
\eeq
It is straightforward to deduce the following lemma from (\ref{e:q1}, \ref{e:M1}, \ref{e:extra}).
\begin{Lemma}\label{l:q1}
Given a pseudo-bilayer $\u_p(x,\veps; {\bf B}_0,\Gamma_0)$ with the mean curvature of the initial interface $\Gamma_0$ far away from zero, the evolution of the pseudo-bilayer in the time scale $\tau_1=\veps t$ satisfies the leading-order evolution system \eqref{e:q1} with
\[
\lim_{\tau_1\to\infty}{\bf B}_1(\tau_1)\cdot{\bf M}=0, \quad \lim_{\tau_1\to\infty}V_{n,0}=0,
\]
where $\lim_{\tau_1\to\infty}{\bf B}_1(\tau_1)=0$ if and only if the leading order term of the mass vector $\mathcal{M}_1$, defined in (\ref{e:calM-def}), is proportional to ${\bf M}$. 
Moreover, any quasi-bilayer profile $\u_q$ associated to $\{{\bf m}_0,\Gamma_0\}$ is an equilibrium of the leading-order evolution system \eqref{e:q1} on time scales up to $\tau_1.$ 
%\[
%{\bf B}_1(\tau_1)\equiv 0, \quad V_{n,0}(s,\tau_1)\equiv 0.
%\]
\end{Lemma}

\begin{Remark}
The convergence of the quantity $E(\tau_1):={\bf B}_1(\tau_1)\cdot {\bf M}$ towards zero has a straightforward physical interpretation: 
for the single-component case, the growing length of the bilayer depletes the amphiphiles in the far-field, or background, region. 
The multicomponent bilayer reaches its maximal length when the inner product of the bilayer mass vector, ${\bf M}$, with the
far field amphiphile density ${\bf B}$ is zero. However it requires
tuning of the initial mass fractions to have all species reach zero in the far-field simultaneously.  That is, while
$\lim_{\tau_1\to \infty}{\bf B}_1(\tau_1)\cdot {\bf M}=0$, it is non-generic that this coincides with the 'emptying' of the far-field region, so
that generically $\lim_{\tau_1\to\infty}{\bf B}_1(\tau_1)\neq 0$. For a two species blend, that is $N=2$, the existence of a one-parameter family of homoclinics, 
as in the universal Birkhoff-billard Example\, \ref{e:2}, would remove this condition as the bilayer composition could adapt continuously to the the available amphiphile supply.
\end{Remark}

%%%%%%%%%%%%%%%%%%%%%%
\paragraph{Time scale $\tau_2=\veps^2 t$: surface-area-preserving Willmore flow}
%%%%%%%%%%%%%%%%%%%%%%
We study the evolution of the set of quasi-bilayers $\mathscr{M}_q$ in the slow time scale $\tau_2=\veps^2t$, assuming that the $O(\veps)$ components of the
background profile ${\bf B}$ has converged to zero on the fast, $\tau_1$, time scale. For the slow, $\tau_2$ time-scale, the initial profile is a quasi-bilayer $\u_q$ associated with interface $\{{\bf m}_0,\Gamma_0\}$. This corresponds
to the inner and outer expansions
\begin{equation}
\tilde{\u}_0(z,s,0)=\phi_h(z), \quad \tilde{\u}_1(z,s,0)=\phi_{h,1}(z), \quad \u_0(x,0)=\u_1(x,0)=0.
\end{equation}
A straightforward but lengthy calculation shows that at the leading-order the solution $u$ of the initial value problem \eqref{e:intlh-1} remains a quasi-bilayer, parmeterized by its
$O(\veps^2)$ back-ground state ${\bf B}_2$ and the interface $\Gamma$, whose evolution is given by, for $\tau_2\geq 0$,
\begin{subnumcases}{\label{e:veps-2} }
\tilde{\u}_0(z,s,\tau_2)\equiv\phi_h(z), \quad \tilde{\u}_1(z,s,\tau_2)\equiv \phi_{h,1}(z), \quad \u_0(x,\tau_2)=\u_1(x,\tau_2)\equiv 0, &\label{e:veps-20}\\
\mu_2(x,\tau_2)=\tilde{\mu}_2(x,\tau_2)=\mathbf{B}_2(\tau_2), &\label{e:veps-21}\\
V_{n,0}(s,\tau_2)=\frac{M_1}{M_2}\left[ \Delta_sH_0+(\frac{\eta_1+\eta_2}{2}+\frac{\mathbf{B}_2\cdot\mathbf{M}}{M_1})H_0-\frac{1}{2}H_0(H_0-a_0)^2-H_1(H_0-a_0) \right], &\label{e:veps-22}
\end{subnumcases}
where we recall $M_1:=\int_\R |\partial_z\phi_h|^2\rmd z$. 
We calculate the $\veps$-scaled total mass vector $\mathcal{M}$ of the amphiphiles, that is,
\[
\mathcal{M}= |\Gamma| \mathbf{M} +\veps\left[ |\Omega |(\nabla_\u^2W(0))^{-2}\mathbf{B}_2+|\Gamma|\int_\R\phi_{h,1}\rmd z+\int_S H_0(s)\rmd s\int_\R z\phi_h(z)\rmd z\right] +\caO(\veps^2).
\]
Given the area expansion $ |\Gamma |(\tau_2)=\gamma_0(\tau_2)+\caO(\veps)$ and the mass expansion $\mathcal{M}=\mathcal{M}_1+\caO(\veps)$, 
it is straightforward to see that 
\beq\label{e:area}
\gamma_0(\tau_2)\equiv |\Gamma_0|.
\eeq
Moreover, plugging the expression of $V_{n,0}$ \eqref{e:veps-22} into the surface area evolution equation \eqref{e:surfevo}, yields the evolution of the background state
\[
\mathbf{B}_2\cdot \mathbf{M}=-M_1\frac{\int_S\left[ -|\nabla_s H_0|^2 +\frac{\eta_1+\eta_2}{2} H_0^2-\frac12 H_0^2(H_0-a_0)^2-H_1H_0(H_0-a_0) \right] \rmd s}{\int_S H_0^2\rmd s}.
\]
Introducing the zero-mass projection operator
\beq\label{e:PiGamma}
\Pi_\Gamma(f)=f-H_0\frac{\int_S fH_0\rmd s}{\int_S H_0^2\rmd s},
\eeq
We conclude that, under the mass conservation, the normal velocity of the interface takes the form
\beq\label{e:nv2}
V_{n,0}=\frac{M_1}{M_2}\Pi_\Gamma\left[  \left(\Delta_s-\frac{1}{2}H_0(H_0-a_0)-H_1  \right)  (H_0-a_0)\right],
\eeq
which, together with the area preserving equation \eqref{e:area}, completes the derivation of Formal Result \ref{t:willmore}.

%%%%%%%%%%%%%%%%%%%%%%%%%%%%%%%%%%%%%%%%%%%%%
\subsection{Evolution of radial bilayers under the mass preserving Willmore flow}
%%%%%%%%%%%%%%%%%%%%%%%%%%%%%%%%%%%%%%%%%%%%%
\begin{figure}[!ht]
  \begin{center}
    \includegraphics[width=1\textwidth]{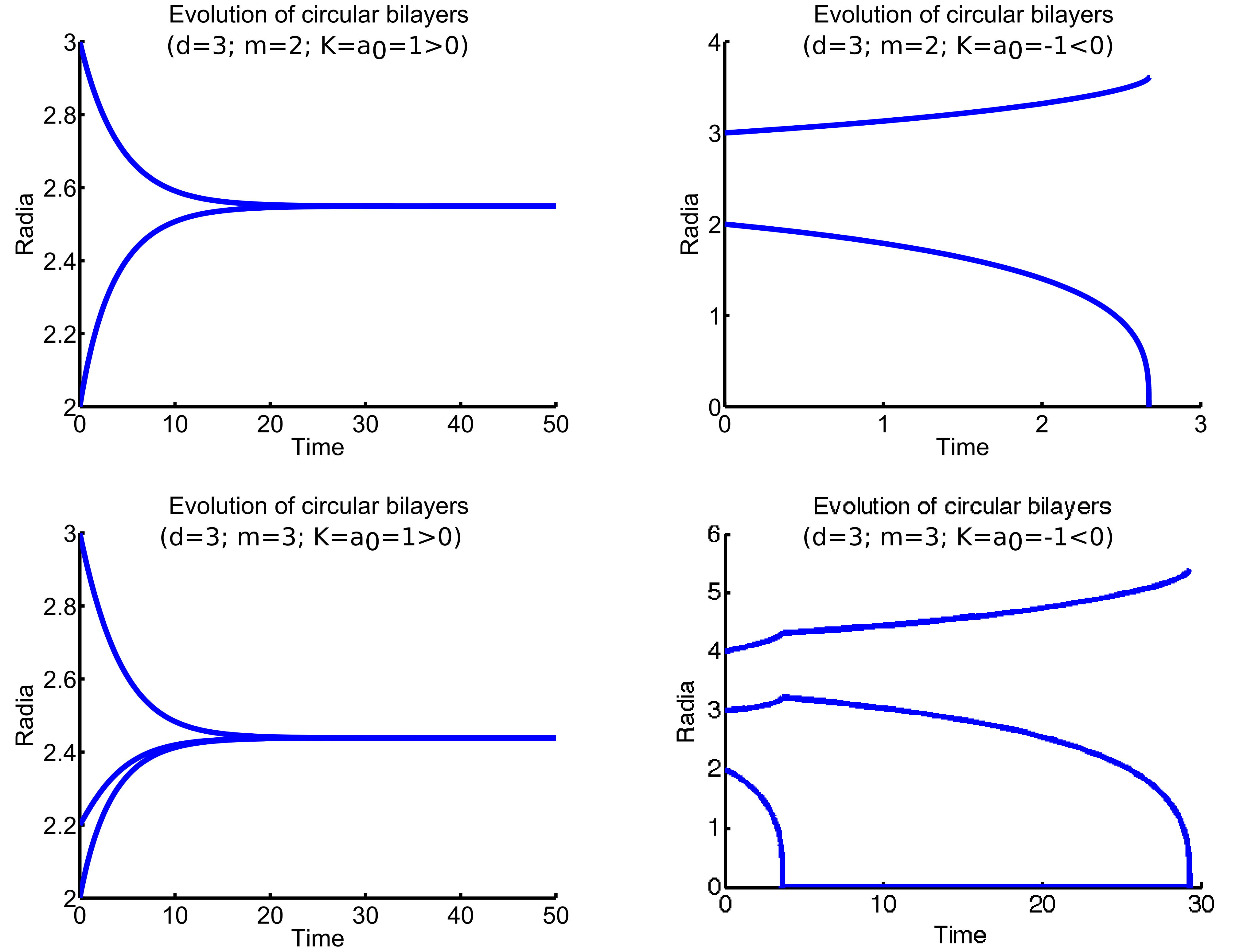}
  \end{center}
  \caption{Simulation of the radius of $m$ spheres with respect to time in spatial dimension $d=3$. For $a_0>0$ the equi-radius
  solution is asymptotically stable, while for $a_0<0$ the smallest radius liposome extinguishes in finite time, leading to an eventual winner-take-all
  scenario.}
\label{f:counterexample}
\end{figure}
To demonstrate the dynamics of the area-preserving Willmore flow, (\ref{e:nv2}), we consider the evolution of 
$m$ well-separated spherically symmetric vesicle cells with radii  $\{R_i(\tau_2)\}_{i=1}^m$ in the domain $\Omega\in\R^d$. For the $i$-th sphere we have 
\[
H_0=(d-1)R_i^{-1}, \quad H_1=-(d-1)R_i^{-2},
\]
and the normal velocity equation \eqref{e:nv2} reduces to 
\beq
\label{e:Nball}
\frac{\rmd R_i}{\rmd \tau_2}=\frac{(d-1)M_1}{2M_2}R_i^{-1}\left\{[(3-d)R_i^{-1}+a_0][(d-1)R_i^{-1}-a_0]-R_c^{-2}\right\},
\eeq
where 
\[R_c^{-2}:=\frac{\sum_{j=1}^{m}[(3-d)R_j^{-1}+a_0][(d-1)R_j^{-1}-a_0]R_j^{d-3}}{\sum_{j=1}^{m}R_j^{d-3}}.\]
Conservation of the total interface area gives a
first integral of the ODE system \eqref{e:Nball}, that is,
\[
\sum_{j=1}^m R_j^{d-1}(\tau_2)=m\overline{R}^{d-1},
\]
where 
\[ \overline{R}:=\sqrt[d-1]{\frac{1}{m}\sum_{j=1}^mR_j^{d-1}(0)}.\] 
Consequently the only equilibrium of the system (\ref{e:Nball}) is the equi-raduis solution $R_i=\overline{R}$ for all $i=1, 2, \cdots, m.$
Linearizing the vector field of the ODE system about the equi-raduis equilbrium we obtain the linearized matrix $\mathbf{L}$ in the form,
\[
\mathbf{L}=-\frac{(d-1)M_1}{M_2 R^{3}}\left[\frac{(d-1)(3-d)}{R}+(d-2)a_0\right](\mathbf{I}_m-\frac{1}{m}\mathbf{J}_m),
\]
where $\mathbf{J}_m$ denotes the $m\times m$ matrix of all ones. The equi-raduis equilibrium is stable for initial data for which 
\beq
\label{e:stableboy}
K(R,a_0,d):=\frac{(d-1)(3-d)}{\overline{R}}+(d-2)a_0>0,
\eeq
and unstable when $K<0$.

In space dimension $d=2$, the intrinsic curvature term $a_0$ drops out of the evolution of (\ref{e:Nball}), which reduces to
\[
\frac{\rmd R_i}{\rmd \tau_2}=\frac{M_1}{2M_2}R_i^{-1}\left(R_i^{-2}-\frac{\sum_{j=1}^{m}R_j^{-3}}{\sum_{j=1}^{m}R_j^{-1}}\right),
\] 
and admits a globally attracting  equi-radius equilibrium $R_i=\overline{R}$, that is, 
\[
\lim_{\tau_2\rightarrow\infty}R_i(\tau)=\overline{R}, \quad i=1,2,\cdots, m.
\]
For $d>2$, the stability of the equi-radius equilibrium depends upon the value of the intrinsic curvature $a_0$. 
More specifically, for $d=3$, we have
\[
\frac{\rmd R_i}{\rmd \tau_2}=2a_0\frac{M_1}{M_2}R_i^{-1}\left(R_i^{-1}-\frac{1}{m}\sum_{j=1}^mR_j^{-1}\right),
\]
where the stability of the equi-radius equilibrium depends uniquely upon the sign of $a_0$. In particular we find 
\[
\lim_{\tau\rightarrow\infty}R_i(\tau_2)=
\begin{cases}
\overline{R}, &a_0> 0,\\
R_i(0), & a_0=0,
\end{cases}
\]
while for $a_0<0$,  the equi-radius equilibrium is unstable. Figure\,\ref{f:counterexample} includes a simulation of the case $a_0<0$ and $m=2$ and $m=3$ 
in  which the smallest circular bilayer shrinks and disappears in finite time, leading to a coarsening phenomenon and eventually a winner-take-all scenario.
%For $d>3$,  the equi-raduisi equilibrium is stable under the evolution (\ref{e:Nball})
%%\[
%%\frac{\rmd R_i}{\rmd \tau}=-\frac{M_1}{2M_2}(d-1)^2(d-3)R_i^{-1}\left\{[R_i^{-1}-\frac{d-2}{(d-1)(d-3)}a_0]^2-\frac{\sum_{j=1}^{m}[R_j^{-1}-\frac{d-2}{(d-1)(d-3)}a_0]^2R_j^{d-3}}{\sum_{j=1}^{m}R_j^{d-3}}
%%\right\},
%%\]
%when 
%\[
%\overline{R}>\frac{(d-3)(d-1)}{(d-2)a_0},
%\]
%and unstable when 
%\[
%\overline{R}<\frac{(d-3)(d-1)}{(d-2)a_0}.
%\]

%%%%%%%%%%%%%%%%%%%%%%%%%%%%%%%%%
\section{Conclusion}
%%%%%%%%%%%%%%%%%%%%%%%%%%%%%%%%%
The multi-component functionalized Cahn-Hilliard (mFCH) free energy is a continuous model for multi-component amphiphilic blends, including plasma membranes, and its gradient
flow describes their slow evolution. We established the existence of quasi-bilayers, a class of co-dimensional one morphologies which encompasses the orbits traced by 
the slow dynamics of bilayer dressings of admissible co-dimension one interfaces, and contains the associated stationary bilayer solutions. We showed that when 
evaluated on the set of quasi-bilayers, mFCH free energy reduces at leading order to the well-known Canham-Helfrich sharp interface free energy, with intrinsic 
curvature determined via a Melnikov parameter arising from the $\veps$-order solenoidal perturbation ${\bf V}$ of the mixing potential ${\bf W}$, (\ref{e:non-solenoidal}). 
For the special class of regularized Birkhoff-billiard mixing potentials, we established the existence of $n$-striation quasi-bilayers for sufficiently small regularization 
parameter $\delta\ll1$, and analyzed their layer-by-layer pearling bifurcation, showing that all order $\caO(\delta^{-2})$ eigenvalues are in one-to-one affiliation with a layer-collision 
and that the associated eigenmodes are localized near the collision point. This analysis establishes the layer-by-layer pearling observed experimentally. Finally we used 
formal multi-scale asymptotics to derive the evolution of quasi-bilayers under the $H^{-1}$ gradient flow of the mFCH,  establishing that on the slow time scale $\tau_2=\veps^2t$, 
the evolution of admissible interfaces follows a mass-preserving Willmore flow. Moreover for radially symmetric bilayers, the intrinsic curvature selects between regimes in which the bilayers seek a common equilibrium radius and a winner-take-all regime in which the smallest radius bilayer is extinguished in finite time. 

Beyond bilayer morphologies, the mFCH free energy accommodates a cornucopia of higher codimensional network structures whose interactions and rich dynamics are unexplored;
moreover the possible complexity of purely bilayer evolution is only hinted at in this work. Indeed, 
%under the regularized Birkhoff-billiard potential, the spectral properties of the  linearized operator around an $n$-striation homoclinic restricted to the 
%nonnegative axis with order smaller than $\caO(\veps^{-2})$ remains to be examined in future. 
of particular interest are Birkhoff-billiard potentials which admit a family of homoclinics parametrized by one or more independent parameters, such as the universal billiard potential
given in Example 3.4; the  associated linearization, (\ref{e:lin-op}) of the underlying dynamical system has a non-trivial kernel, and the associated slow dynamics will couple
the geometric evolution arising from the translational eigenmode to the compositional evolution associated to the hidden symmetries. Taking the functionalization parameter 
$\eta_1$ and $\eta_2$ to depend on $\u$ encodes amphiphile preference for curvature and co-dimension, and yields a dynamic competition in which bilayer composition 
interacts with the geometric evolution and could for example trigger bifurcations leading to endocytosis. 

%It is also important to notice that the mFCH free energy in the current work is a  higher-order perturbation of the scalar FCH free energy, in the sense that we set the ${\bf G}_1$ term in the mGS free energy \eqref{e:mGS},
%\[
%\cF_{mGS}(\u) = \int_\Omega \frac{\veps^4}{2} |{\bf D} \Delta \u|^2 + \veps^2 {\bf G}_1(\u)\cdot {\bf D} \Delta \u + G_2(\u) \rmd x,
%\]
% to be a gradient of a scalar function in the leading order. 
In biological settings, plasma membranes are generically very stable, in particular they are not generically susceptible to pearling bifurcations.
This strict stability to pearling is also generic within the mFCH: indeed the term ${\bf G}_1$ in (\ref{e:mGS}), is not typically the gradient of a scalar 
function, and one anticipates  a more genetical form for the mFCH,
 \beq
\label{e:gmFCH}
\cF_{gM}(\u)=\int_\Omega \frac{1}{2}|\veps^2{\bf D}\Delta \mathbf{u}-{\bf G}(\u)|^2-\veps^2 P(\u)\rmd x,
\eeq
where $\nabla\times {\bf G}\neq 0$.  Assuming the existence of an orbit $\phi_h$ of 
\[
{\bf D}\partial_z^2\u-{\bf G }(\u)=0,
\]
which is homoclinic to zero, the associated linearization $\caL_g:={\bf D}\partial_z^2-\nabla_\u {\bf G}(\phi_h)$ is generically not self-adjoint. 
The second functional derivative of $\cF_{gM}(\u_h)$, evaluated at a dressing, $\u_h$, of an admissible interface with the homoclinic $\phi_h$, admits the expansion,
\[
\mathbb{L}_g:=\frac{\delta^2 \cF_{gM} }{\delta \u^2}(\u_h)=(\caL_g^T+\veps \Delta_s)(\caL_g+\veps \Delta_s)+\caO(\veps).
\]
In many situations \cite{doelman_1998, doelman_2001, doelman_2002}
%Papers 2 and 3  in https://scholar.google.com/citations?user=kB-otT4AAAAJ&hl=en
the spectra of  $\caL_g$  is comprised of eigenvalues $\{\lambda_k\}_{k=1}^N$ with positive real part but with nonzero, order $\caO(1)$ imaginary part.
The corresponding pearling eigenvalues of $\mathbb{L}_g$ take the form $\Lambda_{k,n} = |\lambda_k-\veps^2\beta_n|^2 +\caO(\veps)$, where the real Laplace-Beltrami
eigenvalues $\beta_n>0$ are always an $\caO(1)$ distance away from the strictly complex $\lambda_k$. The lack of positive, purely-real eigenvalues of
$\caL_g$ robustly inhibits pearling and suggests the provocative statement: the choices available in the packing of lipids within multicomponent membranes generically 
serves to inhibit the mechanisms of pearling. This simple analysis is in sympathy with diverse observations of plasma membranes; and suggests that dynamical systems 
has a potentially significant role to play in our understanding of the nature of these essential components of cellular design.

%%%%%%%%%%%%%%%%%%%%%%%%%%%%%%%%%%%%%%%%%%%%%%%%%%%%%%
% Literature
%%%%%%%%%%%%%%%%%%%%%%%%%%%%%%%%%%%%%%%%%%%%%%%%%%%%%%
%\bibliographystyle{siam}

\end{document}